\documentclass{icmart}
\usepackage[all]{xy}

\contact[andrus@famaf.unc.edu.ar]{FaMAF,
Universidad Nacional de C\'ordoba. CIEM -- CONICET. \newline \noindent
Medina Allende s/n (5000) Ciudad Universitaria, C\'ordoba,
Argentina}

\numberwithin{equation}{section}

\newtheorem{theorem}{Theorem}[section]

\newtheorem{conjecture}[theorem]{Conjecture}

\theoremstyle{definition}
\newtheorem{definition}[theorem]{Definition}
\newtheorem{example}[theorem]{Example}

\newcommand\ad{\operatorname{ad}}
\newcommand{\Aut}{\operatorname{Aut}}
\newcommand{\gr}{\operatorname{gr}}

\newcommand{\id}{\operatorname{id}}
\newcommand\im{\operatorname{im}}
\newcommand{\Ind}{\operatorname{Ind}}

\newcommand{\Irr}{\operatorname{Irr}}

\newcommand{\ord}{\operatorname{ord}}

\newcommand\re{\operatorname{re}}
\newcommand\rk{\operatorname{rk}}
\newcommand\sgn{\operatorname{sgn}}
\newcommand\supp{\operatorname{supp}}

\newcommand{\trid}{\triangleright}

\newcommand{\toba}{{\mathfrak B}}
\newcommand{\C}{{\mathfrak C}}
\def\J{\mathfrak{J}}
\def\h{\mathfrak{h}}
\def\g{\mathfrak{g}}
  
\newcommand{\ug}{\mathfrak{u}}  
\newcommand{\vg}{\mathfrak{v}}  
\newcommand{\A}{\ug}

\newcommand{\ba}{\mathbf{a} }

\newcommand{\GL}{\mathbf{GL}}
\newcommand\ku{\mathbf{k}}
\newcommand{\bm}{\mathbf{m} }
\newcommand{\n}{\mathbf{n}}
\newcommand{\bp}{\mathbf{p} }
\newcommand{\PSL}{\mathbf{PSL}}
\newcommand{\bq}{\mathbf{q} }
\newcommand{\bM}{\mathbf{M} }

\def\cC{\mathcal{C}}
\newcommand\D{\mathcal D}
\newcommand{\Ee}{{\mathcal E}}
\def\cG{\mathcal{G}}

\newcommand{\Oc}{{\mathcal O}}
\def\Pc{\mathcal{P}}
\def\Q{\mathcal{Q}}
\def\cR{\mathcal{R}}
\newcommand\Ss{\mathcal{S}}
\newcommand\T{\mathcal{T}}
\def\cW{\mathcal{W}}
\def\cX{\mathcal{X}}

\newcommand\am{\mathbb A_m}
\newcommand{\kk}{\mathbb C}
\newcommand\dm{\mathbb D_m}
\newcommand\F{\mathbb F}
\newcommand\G{\mathbb G}
\newcommand\I{\mathbb I}
\newcommand\N{\mathbb N}
\newcommand\Sb{\mathbb S}
\newcommand\st{\mathbb S_3}
\newcommand\sk{\mathbb S_4}
\newcommand\sco{\mathbb S_5}
\newcommand\sei{\mathbb S_6}

\newcommand\sm{\mathbb S_m}
\newcommand{\Z}{\mathbb Z}
\newcommand\zt{\mathbb{Z}^\theta}

\newcommand\varep{\varepsilon}
\newcommand{\VGamma}{\widehat{\Gamma}}  

\newcommand{\ydh}{{}_{H}^{H}\mathcal{YD}}
\newcommand{\ydl}{{}_{L}^{L}\mathcal{YD}}
\newcommand{\ydho}{{}_{H_{[0]}}^{H_{[0]}}\mathcal{YD}}
\newcommand{\ydk}{{}_{K}^{K}\mathcal{YD}}

\newcommand{\ydg}{{}_{\Gamma}^{\Gamma}\mathcal{YD}}
\newcommand{\ydgg}{{}_{G}^{G}\mathcal{YD}}
\newcommand{\mc}[2]{{}_{#1}^{#2}\mathcal{M}}

\newcommand{\xij}[1]{a_{(#1)}}

\title[On finite-dimensional Hopf algebras]{On finite-dimensional Hopf algebras}

\author[Nicol\'as Andruskiewitsch]
{Nicol\'as Andruskiewitsch\thanks{This work was partially supported by
ANPCyT-Foncyt, CONICET, Secyt (UNC).}}

\begin{document}

\begin{center}
 \emph{Dedicado a Biblioco 34}
\end{center}

\begin{abstract}
This is a survey on the state-of-the-art of the classification of finite-dimensio\-nal complex Hopf algebras. 
This general question is addressed through the consideration of different classes of such Hopf algebras.
Pointed Hopf algebras constitute the class best understood; the classification of those with abelian group is expected to be
completed soon and there is substantial progress in the non-abelian case. 

\end{abstract}



\maketitle

\section{Introduction}

Hopf algebras were introduced in the 1950's from three different perspectives: algebraic groups in positive characteristic,
cohomology rings of Lie groups, and group objects in the category of von Neumann algebras. 
The study of non-commutative non-cocommutative Hopf algebras started in the 1960's. 
The fundamental breakthrough is  Drinfeld's report \cite{Dr1}. Among many contributions and ideas,
a systematic construction of solutions of the quantum Yang-Baxter equation (qYBE) was presented. Let $V$ be a vector space.
The qYBE is equivalent to the braid equation: 
\begin{align}\label{eq:braid}
(c\otimes \id) (\id\otimes c) (c\otimes \id) &= (\id\otimes c) (c\otimes \id) 
(\id\otimes c),& 
c&\in GL(V\otimes V).
\end{align}
If $c$ satisfies \eqref{eq:braid}, then $(V,c)$ is called a braided vector space;  this is a down-to-the-earth version
of a braided tensor category \cite{JS}.
Drinfeld introduced the notion of quasi-triangular Hopf algebra, meaning a pair $(H, R)$ where $H$ is a Hopf algebra and 
$R\in H\otimes H$ is invertible and satisfies the approppriate conditions, so that every $H$-module $V$ becomes a braided
vector space, with $c$ given by the action of $R$ composed with the usual flip.
Furthermore, every finite-dimensional Hopf algebra $H$ gives rise to a quasi-triangular Hopf algebra, namely the Drinfeld double $D(H) = H \otimes H^*$
as vector space. If $H$ is not finite-dimensional, some precautions have to be taken to construct $D(H)$, or else one considers 
Yetter-Drinfeld modules, see \S \ref{subsec:repns}. In conclusion, every Hopf algebra is a source of solutions of the braid equation.
Essential examples of quasi-triangular Hopf algebras are the quantum groups $U_q(\mathfrak g)$ \cite{Dr1, J} and the
finite-dimensional variations $\ug_q(\mathfrak g)$  \cite{L1, L2}. 

In the approach to the classification of Hopf algebras exposed in this report, braided vector spaces and braided tensor categories play a decisive role;
and the finite quantum groups are the main actors in one of the classes that splits off.

By space limitations, there is a selection of the topics and references included.
Particularly, we deal with finite-dimensional Hopf algebras over an algebraically closed field of characteristic zero with special emphasis on
description of examples and classifications. 
Interesting results on Hopf algebras either infinite-dimensional, or over other fields, unfortunately can not be reported.
There is no account of the many deep results on tensor categories, see \cite{EGNO}. 
Various basic fundamental results are not explicitly cited, we refer to \cite{bariloche, Ma3, Mo1,  Rad, Sch2, Sw} for them; classifications of Hopf algebras
of fixed dimensions are not evoked, see \cite{BG, Ng, Z}.

\smallbreak\subsubsection*{Acknowledgements} I thank all my coauthors for pleasant and fruitful collaborations. 
I am particularly in debt with Hans-J\"urgen Schneider, Mat\'\i as Gra{\~n}a and Iv\'an Angiono for sharing their 
ideas on pointed Hopf algebras with me.

\section{Preliminaries}

Let $\theta\in\mathbb N$ and $\I= \I_\theta = \{1,2,\dots,\theta\}$. The base field is $\kk$. If $X$ is a set, then
$\vert X\vert$ is its cardinal and $\kk X$ is the vector space with basis $(x_i)_{i\in X}$.
Let $G$ be a group: we denote by $\Irr G$ the set of isomorphism classes of irreducible
representations of $G$ and by $\widehat{G}$ the subset of those of dimension 1;  by $G^x$ the centralizer of $x\in G$; and by $\Oc_x^G$ its conjugacy class.
More generally we denote by $\Irr \C$ the set of isomorphism classes of simple objects in an abelian category $\C$.
The group of $n$-th roots of 1 in $\kk$ is denoted $\G_n$; also $\G_{\infty} = \bigcup_{n\ge1} \G_n$.
The group presented by $(x_i)_{i\in I}$ with relations $(r_j)_{j\in J}$ is denoted $\langle (x_i)_{i\in I} \vert (r_j)_{j\in J} \rangle$.
The notation for Hopf algebras is  standard: $\Delta$, $\varep$, $\Ss$,  denote
respectively the comultiplication, the counit, the antipode (always assumed bijective, what happens in the finite-dimensional case). We use Sweedler's
notation:  $\Delta(x) =  x_{(1)} \otimes x_{(2)}$. 
Similarly, if $C$ is a coalgebra and $V$ is a left comodule with structure map 
$\delta: V \to C \otimes V$, then $\delta(v) = v_{(-1)} \otimes v_{(0)}$.  
If $D, E$ are subspaces  of $C$, then  $D \wedge E = \{c\in C: \Delta(c) \in D\otimes C + C \otimes E\}$; also 
$\wedge^{0} D = D$ and $\wedge^{n+1} D = (\wedge^{n} D) \wedge D$ for $n>0$. 

\smallbreak\subsection{Basic constructions and results}
The first examples of finite-dimensional Hopf algebras are the group algebra $\kk G$ of a finite group $G$
and its dual, the algebra of functions $\kk^G$. Indeed, 
the dual of a finite-dimensional Hopf algebra is again a Hopf algebra by transposing operations.
By analogy with groups, several authors explored the notion of extension of Hopf algebras at various 
levels of generality; in the finite-dimensional context, every extension 
$\kk \to A \rightarrow C \rightarrow B \to \kk$ can be described as $C$ with underlying vector space
$A\otimes B$, via a heavy machinery of actions, 
coactions and non-abelian cocycles, but actual examples are rarely found in this way (extensions from a different perspective 
are in \cite{AG-compo}).  
Relevant exceptions are the so-called \emph{abelian extensions} \cite{K-ext} (rediscovered by Takeuchi and Majid):
here the input is a matched pair of groups $(F, G)$
with mutual actions $\triangleright$, $\triangleleft$ (or equivalently, an exact factorization of a finite group).  
The actions give rise to a Hopf algebra  $\kk^G \# \kk F$. 
The multiplication and comultiplication can be further modified by  compatible cocycles
$(\sigma, \kappa)$, producing to the abelian extension 
$\kk \to \kk^G \rightarrow \kk^G {}^\kappa\hspace{-3pt}\#_\sigma \kk F \rightarrow \kk F \to \kk$.
Here $(\sigma, \kappa)$ turns out to be a 2-cocycle in the total complex associated to a double complex built
from the matched pair; the relevant $H^2$ is computed via the so-called Kac exact sequence.

\smallbreak
It is natural to approach Hopf algebras by considering algebra or coalgebra invariants.
 There is  no preference in the finite-dimensional setting but coalgebras and comodules are locally finite, so we privilege the coalgebra ones
 to lay down general methods. The basic coalgebra invariants of a Hopf algebra $H$ are:

\smallbreak \noindent $\circ$ The group $G(H) =\{g\in H - 0: \Delta(g) = x\otimes g\}$ of group-like elements of $H$.

\smallbreak \noindent $\circ$ The space of skew-primitive elements $\Pc_{g,h}(H)$, $g, h\in G(H)$; $\Pc(H) := \Pc_{1,1}(H)$.

\smallbreak \noindent $\circ$ The coradical $H_0$, that is the sum of all simple subcoalgebras.

\smallbreak \noindent $\circ$ The coradical filtration $H_0 \subset H_1 \subset \dots$, where $ H_n= \wedge^{n} H_0$; then $H = \bigcup_{n \geq 0} H_n$.
 
\smallbreak\subsection{Modules}\label{subsec:repns}
The category $\mc{H}{}$ of left modules over a Hopf algebra $H$ is monoidal with tensor product defined by the comultiplication; 
ditto for the category $\mc{}{H}$ of left comodules, with tensor product defined by the multiplication. Here are two ways 
to deform Hopf algebras without altering one of these categories.

\smallbreak \noindent $\bullet$  Let $F\in H\otimes H$ be invertible such that $(1\otimes F) (\id \otimes \Delta)(F) = 
(F\otimes 1) (\Delta \otimes \id)(F)$ and $(\id \otimes \varepsilon)(F) = 
(\varepsilon \otimes \id)(F) = 1$. Then $H^{F}$ (the same 
algebra with comultiplication $\Delta^{F}:= F \Delta F^{-1}$) is again a Hopf algebra,
named the \emph{twisting} of $H$ by $F$ \cite{Dr-quasi}. The monoidal categories $\mc{H}{}$ and $\mc{H^F}{}$ are equivalent. 
If $H$ and $K$ are finite-dimensional Hopf algebras with $\mc{H}{}$ and $\mc{K}{}$ equivalent as monoidal categories, 
then there exists $F$ with $K\simeq H^F$ as Hopf algebras (Schauenburg, Etingof-Gelaki). 
Examples of twistings not mentioned elsewhere in this report are in \cite{EN, Momb}.

\smallbreak \noindent $\bullet$ Given a linear map $\sigma: H\otimes H \to \kk$ with analogous conditions, there is a Hopf algebra
$H_{\sigma}$ (same coalgebra, multiplication twisted by $\sigma$) such that the monoidal categories $\mc{}{H}$ and $\mc{}{H_\sigma}$ are equivalent \cite{DT}.
 
\smallbreak
A Yetter-Drinfeld module $M$ over $H$ is left $H$-module and left $H$-comodule with the compatibility
$\delta(h.m) = h_{(1)} m_{(-1)}\Ss(h_{(3)}) \otimes h_{(2)} \cdot  m_{(0)}$, for all $m\in M$ and $h\in H$.
The category $\ydh$ of Yetter-Drinfeld modules is braided monoidal.\label{page:yd} 
That is, for every $M,N \in \ydh$,  there is a natural isomorphism $c: M\otimes N \to N\otimes 
M$ given by $c(m\otimes n) = m_{(-1)}\cdot n \otimes m_{(0)}$, $m\in M$, $n\in N$.
When $H$ is finite-dimensional, the category $\ydh$ is equivalent, as a braided monoidal category, to  $\mc{D(H)}{}$.

\smallbreak
The definition of Hopf algebra makes sense in any braided monoidal category. Hopf algebras in $\ydh$ are interesting because of the following 
facts--discovered by Radford and interpreted categorically by Majid, see \cite{Ma3, Rad}:

\smallbreak \noindent $\diamond$ If $R$ is a Hopf algebra in $\ydh$, then $R\# H := R \otimes H$ with semidirect product and coproduct is a Hopf algebra,
named the \emph{bosonization} of $R$ by $H$.

\smallbreak \noindent $\diamond$ Let $\pi, \iota$ be Hopf algebra maps as in
$\xymatrix{K \ar@{->>}_{\pi}[r] & \ar@/^-1pc/_{\iota}[l] H}$ with $\pi\iota = \id_H$. 
Then $R = H^{\text{co} \pi} := \{x\in K: (\id \otimes \pi)\Delta(x) = 1\otimes x\}$ is a Hopf algebra in $\ydh$ and $K \simeq R\# H$.

\smallbreak 
For instance, if $V\in \ydh$, then the tensor algebra $T(V)$ is a Hopf algebra in $\ydh$, by requiring 
$V \hookrightarrow \Pc(T(V))$. If $c:V\otimes V \to V\otimes V$ satisfies $c = - \tau$,  $\tau$ the usual flip, then the exterior algebra $\Lambda(V)$ 
is a Hopf algebra in $\ydh$.

There is a  braided adjoint action of a Hopf algebra $R$ in $\ydh$ on itself, see e.g. \cite[(1.26)]{AHS}. If $x\in \Pc(R)$ and $y \in R$, 
then $\ad_c(x)(y) =  xy - \text{mult } c(x\otimes y)$.

\smallbreak\subsubsection{Triangular Hopf algebras}\label{subsec:triangular}
A quasitriangular Hopf algebra $(H,R)$ is \emph{triangular} if the braiding induced by $R$ is a symmetry:
$c_{V\otimes W}c_{W\otimes V} = \id_{W\otimes V}$ for all $V, W \in \mc{H}{}$. A finite-dimensional
triangular Hopf algebra is a twisting of a bosonization $\Lambda(V)\# \kk G$, where $G$ is a finite group and $V\in \ydgg$ has $c = -\tau$ \cite{AEG}.
This lead eventually to the classification of triangular finite-dimensional  Hopf algebras \cite{EG-clas-triang}; previous work on the semisimple
case culminated in \cite{EG-clas-triang-ss}.

\smallbreak\subsection{Semisimple Hopf algebras}

The algebra of functions $\kk^G$ on a finite group $G$ admits a Haar measure, 
i.e., a linear function $\smallint: \kk^G \to \kk$ 
invariant under left and right translations, namely $\smallint = $ 
sum of all elements in the standard basis of $\kk G$. 
This is adapted as follows: a right integral \emph{on} a Hopf algebra $H$ is 
a linear function $\smallint: H \to \kk$ 
which is invariant under the left  regular coaction: analogously there is the notion of left integral. 
The notion has various applications. Assume that $H$ is finite-dimensional.
Then an integral \emph{in} $H$ is an integral on $H^*$; the subspace of left integrals in $H$ has dimension one, and there is  a generalization of
 Maschke's theorem for finite groups: $H$ is semisimple if and only if $\varep(\Lambda)\neq 0$ for any  integral $0\neq \Lambda\in H$.
This characterization of semisimple Hopf algebras, valid in any characteristic, is one of several, some valid only in characteristic 0.
See \cite{Sch2}. Semisimple Hopf algebras can be obtained as follows:

\smallbreak \noindent $\diamond$ A finite-dimensional Hopf algebra $H$ is semisimple if and only if it is cosemisimple (that is, $H^*$ is semisimple).

\smallbreak \noindent $\diamond$ Given an extension $\kk \to K \rightarrow H \rightarrow L \to \kk$, $H$ is semisimple iff $K$ and $L$ are.
Notice that there are semisimple extensions that are not abelian \cite{GNN, Nat-ext, Nik}.

\smallbreak \noindent $\diamond$ If $H$ is semisimple, then so are $H^{F}$ and $H_{\sigma}$, for any twist $F$ and cocycle $\sigma$.
If $G$ is a finite simple group, then any twisting of $\kk G$ is a simple Hopf algebra (i.e., not a non-trivial extension) \cite{Nik-tw}, 
but the converse is not true \cite{GN}.

\smallbreak \noindent $\diamond$ A bosonization $R\# H$ is semisimple iff $R$ and $H$ are.

To my knowledge, all examples of semisimple Hopf algebras arise from group algebras by the preceding constructions; this was proved in 
\cite{N4-mem, N4-60} in low dimensions and in \cite{ENO} for dimensions $p^aq^b$, $pqr$, where $p$, $q$ and $r$ are primes. 
See \cite[Question 2.6]{bariloche}. An analogous question in terms of fusion categories: is any semisimple Hopf algebra weakly group-theoretical? See
\cite[Question 2]{ENO}.

There are only finitely many isomorphism classes of \emph{semisimple} Hopf algebras in each dimension \cite{St-finite-ss}, but this fails in general
\cite{AS-jalg, BDG}.

\begin{conjecture} (Kaplansky). Let $H$ be a semisimple Hopf algebra. The dimension of every $V\in \Irr \mc{H}{}$ divides the dimension of $H$.
 \end{conjecture}

The answer is affirmative for iterated extensions of group algebras and duals of group algebras \cite{MW} 
and notably for semisimple quasitriangular Hopf algebras \cite{EG-frob-prop}.

\section{Lifting methods}\label{sec:methods}

\smallbreak\subsection{Nichols algebras}\label{subsec:nichols}

Nichols algebras are a special kind of Hopf algebras in braided tensor categories. 
We are mainly interested in Nichols algebras in the braided category $\ydh$, where $H$ is a
Hopf algebra, see page \pageref{page:yd}. In fact, there is a functor $V\mapsto \toba(V)$ from $\ydh$ to the category of 
Hopf algebras in $\ydh$. Their first appearence is in the precursor  \cite{N}; they were
rediscovered in   \cite{W} as part of a ``quantum differential calculus", and in \cite{Lu}
to present the positive part of $U_q(\mathfrak g)$. 
See also \cite{Ro2, Sbg}. 

There are several, unrelated at the first glance, alternative definitions.
Let $V\in \ydh$. The first definition uses the representation of the braid group $\mathbb B_n$ in $n$ strands on $V^{\otimes n}$,
given by $\varsigma_i \mapsto \id\otimes c\otimes \id$, $c$ in $(i, i+1)$ tensorands; here recall that
$\mathbb B_n = \langle\varsigma_1, \dots, \varsigma_{n-1}\vert \varsigma_i\varsigma_j = \varsigma_j\varsigma_i, \vert i-j\vert >1, 
\varsigma_i\varsigma_j\varsigma_i = \varsigma_j\varsigma_i\varsigma_j, \vert i-j\vert = 1\rangle$. Let 
$M: \mathbb S_n \to\mathbb B_n$ be the Matsumoto section and let $\Q_n :V^{\otimes n} \to V^{\otimes n}$ be the
quantum symmetrizer, $\Q_n = \sum_{s\in \mathbb S_n} M(s) :V^{\otimes n} \to V^{\otimes n}$. Then define
\begin{align}\label{eq:nichols-first-def}
\J^n(V) &= \ker \Q_n, & \J(V) &= \oplus_{n\geq 2} \J^n(V),& \toba(V) &= T(V) / \J(V).
\end{align}
Hence $\toba(V)= \oplus_{n\geq 0}\toba^n(V)$ is a graded Hopf algebra in $\ydh$ with $\toba^0(V) = \kk$,  $\toba^1(V) \simeq V$;
by \eqref{eq:nichols-first-def} the algebra structure depends only on $c$. To explain the second definition, let us observe that the tensor
algebra $T(V)$ is a Hopf algebra in $\ydh$ with comultiplication determined by $\Delta(v) = v\otimes 1 + 1 \otimes v$ for $v\in V$.
Then $\J(V)$ coincides with the largest homogeneous ideal of $T(V)$  generated by elements of degree $\ge 2$ that is also a coideal.
Let now 
$T= \oplus_{n\geq 0}T^n$ be a graded Hopf algebra in $\ydh$ with $T^0 = \kk$. Consider the conditions
\begin{align}\label{eq:nichols-charac-alg}
   T^1 &  \text{ generates $T$ as an algebra},  \\ \label{eq:nichols-charac-coalg}
   T^1 &=  \Pc (T).
\end{align}
These requirements are dual to each other: if $T$ has finite-dimensional homogeneous components and $R = \oplus_{n\geq 0}R^n$ is the graded
dual of $T$, i.e., $R^n = (T^n)^*$, then $T$ satisfies \eqref{eq:nichols-charac-alg} if and only if 
$R$ satisfies \eqref{eq:nichols-charac-coalg}. These conditions determine $\toba(V)$ up to isomorphisms, 
as the unique graded connected Hopf algebra $T$ in $\ydh$ that satisfies  $T^1 \simeq V$, \eqref{eq:nichols-charac-alg} and \eqref{eq:nichols-charac-coalg}.
There are still other characterizations of $\J(V)$, e.g.  as the radical of a suitable homogeneous bilinear form on $T(V)$,
or as the common kernel of some suitable skew-derivations.  See \cite{AS-cambr} for more details.

Despite all these different definitions, Nichols algebras are extremely difficult to deal with, e.g. to present by generators and relations,
or to determine when a Nichols algebra has finite dimension or finite Gelfand-Kirillov dimension. 
It is not even known a priori whether the ideal $\J(V)$ is finitely generated, except in a few specific cases. 
For instance, if $c$ is a symmetry, that is $c^2 =\id$, or satisfies a Hecke condition with generic parameter, then $\toba(V)$ is quadratic.
By the efforts of various authors, we have some understanding 
of finite-dimensional Nichols algebras of braided vector spaces either of diagonal or of rack type, see \S \ref{subsec:nichols-diagonal},
\ref{subsec:nichols-rack}.

\smallbreak\subsection{Hopf algebras with the (dual) Chevalley property}\label{subsec:chevalley}
We now explain how Nichols algebras enter into our approach to the classification of Hopf algebras. Recall that
a Hopf algebra has the \emph{dual  Chevalley property} if the tensor product of two simple comodules is semisimple, or 
equivalently if its coradical is a (cosemisimple) Hopf subalgebra. For instance, a \emph{pointed} Hopf algebra,
one whose simple comodules have all dimension one, has the dual Chevalley property and its coradical is a group algebra.
Also, a \emph{copointed} Hopf algebra (one whose coradical is the algebra of functions on a finite group)
has the dual Chevalley property.
The Lifting Method is formulated in this context \cite{AS-jalg}. Let $H$ be a Hopf
algebra with the dual Chevalley property and set $K := H_0$. Under this assumption, the graded coalgebra 
$\gr H = \oplus_{n\in \mathbb N_0}\gr^n H$ associated to the coradical filtration becomes a Hopf algebra and 
considering the homogeneous projection $\pi$ as in
$\xymatrix{R = H^{\text{co}\, \pi} \ar@{^{(}->}[r] & \gr H  \ar@<-0.5ex>@{->}_{\quad  \pi}[r] & \ar@<-0.5ex>[l] K}
$
we see that $\gr H \simeq R \# K$. The subalgebra of coinvariants $R$ is a graded Hopf algebra in $\ydk$
that inherits the grading with $R^0 = \kk$; it satisfies \eqref{eq:nichols-charac-coalg} since the grading comes from the coradical filtration.
Let $R'$ be the subalgebra of $R$ generated by  $ V :=R^1$; then $R' \simeq \toba (V)$. The braided vector space $V$ is a basic invariant
of $H$ called its \emph{infinitesimal braiding}.
Let us fix then a semisimple Hopf algebra $K$.
To classify all finite-dimensional Hopf algebras $H$ with $H_0 \simeq K$ as Hopf algebras,  we have to address the 
following questions.

\begin{enumerate}\renewcommand{\theenumi}{(\alph{enumi})}\itemsep1pt \parskip0pt \parsep0pt
\item\label{item:method-nichols}  Determine those $V \in \ydk$ such that $\mathfrak  B(V)$  is finite-dimensional, 
and give an efficient defining set of relations of these.

\item\label{item:method-gen-deg-one}  Investigate whether any {\it finite-dimensional} graded  Hopf 
algebra $R$ in $\ydk$ 
satisfying $R^0 = \kk$ and $P(R) = R^1$,  is a Nichols algebra.

\item\label{item:method-lifting}  Compute all Hopf algebras $H$ such that $\gr H \simeq \mathfrak  B(V) \# K$, $V$ as in \ref{item:method-nichols}.
\end{enumerate}

Since the Nichols algebra $\toba(V)$ depends as an algebra (and as a  coalgebra) only on the braiding $c$,
it is convenient to restate Question \ref{item:method-nichols} as follows:

\begin{enumerate}\itemsep1pt \parskip0pt \parsep0pt
\item[(a$_1$).] Determine those braided vector spaces $(V, c)$ in a suitable class
such that $\dim \mathfrak  B(V) < \infty$,  and give an efficient defining set of relations of these.

\item[(a$_2$).] For those $V$ as in (a$_1$), find in how many ways, if any, they can be  realized as Yetter-Drinfeld modules over $K$.
\end{enumerate}

For instance, if $K = \kk\Gamma$,  $\Gamma$ a finite abelian group, then the suitable class is that of braided vector spaces of diagonal type.
In this context, Question (a$_2$) amounts to solve systems of equations in $\Gamma$.
The answer to  \ref{item:method-nichols} is instrumental to attack \ref{item:method-gen-deg-one} and \ref{item:method-lifting}.
Question \ref{item:method-gen-deg-one}  can be rephrased in two equivalent statements:

\begin{enumerate}\itemsep1pt \parskip0pt \parsep0pt
\item[(b$_1$).] Investigate whether any {\it finite-dimensional} graded  Hopf 
algebra $T$ in $\ydk$  with $T^0 = \kk$ and generated as algebra by $T^1$,  is a Nichols algebra.

\item[(b$_2$).] Investigate whether any {\it finite-dimensional}  Hopf 
algebra $H$  with $H_0 = K$ is generated as algebra by $H_1$.
\end{enumerate}

We believe that the answer to \ref{item:method-gen-deg-one}  is affirmative at least when $K$ is a group algebra.
In other words, by the reformulation (b$_2$):

\begin{conjecture}\label{conj:gendegone} \emph{\cite{AS-adv}}
 Every {\it finite-dimensional} pointed Hopf algebra is generated by group-like and skew-primitive elements.
\end{conjecture}

As we shall see in \S \ref{subsec:gen-deg-1}, 
the complete answer to \ref{item:method-nichols} is needed in the approach proposed in \cite{AS-adv} to attack Conjecture \ref{conj:gendegone}.
It is plausible that the answer of (b$_2$) is affirmative for every semisimple Hopf algebra $K$.
Question \ref{item:method-lifting}, known as \emph{lifting} of the relations, also
requires the knowledge of the generators of $\J(V)$, see \S \ref{subsec:deformations}.

\smallbreak\subsection{Generalized Lifting method}\label{subsec:generalized-lifting-method}

Before starting with the analysis of the various questions in \S \ref{subsec:chevalley}, we discuss a possible approach to more general Hopf algebras \cite{AC}.
Let $H$ be a Hopf algebra; we consider the following invariants of $H$:

\smallbreak \noindent $\circ$ The \emph{Hopf coradical} $H_{[0]}$ is the subalgebra generated by $H_0$.  

\smallbreak \noindent $\circ$ The \emph{standard filtration} $H_{[0]} \subset H_{[1]} \subset \dots$,  $H_{[n]}=\wedge^{n+1} H_{[0]}$; 
then $H = \bigcup_{n \geq 0} H_{[n]}$.

\smallbreak 
If $H$ has the dual Chevalley property, then $H_{[n]} = H_n$ for all $n\in \mathbb N_0$. 
In general,  $H_{[0]}$ is a Hopf subalgebra of $H$ with coradical $H_0$ and we may consider the  graded Hopf 
algebra $\gr H=\oplus_{n \geq 0} H_{[n]}/H_{[n-1]}$. 
As before, if $\pi:\gr H \rightarrow H_{[0]}$ is the homogeneous projection, then  $R = (\gr H)^{\text{co}\, \pi}$ is a 
Hopf algebra in $\ydho$ and $\gr H \cong R \# H_{[0]}$.
Furthermore, $R= \oplus_{n \geq 0} R^n$ with grading inherited from $\gr H$.
This discussion raises the following questions.

\begin{enumerate}\renewcommand{\theenumi}{(\Alph{enumi})}\itemsep1pt \parskip0pt \parsep0pt
\item\label{item:gen-method-coalgebras}  
Let $C$ be a finite-dimensional cosemisimple coalgebra and $S: C\to C$ a bijective 
anti-coalgebra map. 
Classify all finite-dimensional  Hopf algebras $L$ generated by $C$, such that $\Ss {\vert_C} = S$.

\item\label{item:gen-method-R} Given $L$ as in the previous item, classify all finite-dimensional connected graded Hopf algebras $R$ in $\ydl$.

\item\label{item:gen-method-lifting}  Given $L$ and $R$ as in previous items, classify all deformations or liftings, that is, classify all Hopf algebras $H$
such that $\gr H \cong R\# L$.
\end{enumerate}

Question \ref{item:gen-method-coalgebras} is largely open, except for the remarkable \cite[Theorem 1.5]{St}:
if $H$ is a Hopf algebra generated by  an $\Ss$-invariant 4-dimensional simple subcoalgebra $C$, such that  
$1 < \ord (\Ss^{2}\hspace{-1pt}\vert_{_C}) < \infty$, 
then $H$ is a Hopf algebra quotient of the quantized algebra of functions on $SL_2$ at a root of unity $\omega$.
Nichols algebras enter into the picture in Question \ref{item:gen-method-R}; if $V = R^1$, then $\toba(V)$ is a subquotient of $R$.
Question \ref{item:gen-method-lifting} is completely open, as it depends on the previous Questions.

\smallbreak\subsection{Generalized root systems and Weyl groupoids}\label{subsec:grs}

Here we expose two important notions introduced in \cite{HY}.

Let $\theta \in \N$ and $\I = \I_{\theta}$.
A \emph{basic datum} of type $\I$ is a pair  $(\cX, \rho)$, where $\cX \neq \emptyset$ is a  set and
$\rho: \I  \to \Sb_{\cX}$ is a map such that $\rho_i^2 = \id$ for all $i\in \I$.
 Let $\Q_{\rho}$ be the quiver $\{\sigma_i^x := (x, i, \rho_i(x)): i\in \I, x\in \cX\}$ over $\cX$, with
$t(\sigma_i^x) = x$, $s(\sigma_i^x) = \rho_i(x)$ (here $t$ means target, $s$ means source). Let $F(\Q_{\rho})$ be the free groupoid over $\Q_{\rho}$;
in any quotient of $F(\Q_{\rho})$, we denote 
\begin{align}\label{eq:conventio}
\sigma_{i_1}^x\sigma_{i_2}\cdots \sigma_{i_t} =
\sigma_{i_1}^x\sigma_{i_2}^{\rho_{i_1}(x)}\cdots \sigma_{i_t}^{\rho_{i_{t-1}} \cdots \rho_{i_1}(x)};
\end{align}
i.e., the implicit superscripts are those  allowing compositions.

\smallbreak\subsubsection{Coxeter groupoids}\label{subsubsec:cox-grpds}
A \emph{Coxeter datum} is a triple $(\cX, \rho, \bM)$, where $(\cX, \rho)$ is a basic datum of type $\I$
and $\bM = (\bm^x)_{x\in \cX}$  is a family of Coxeter matrices $\bm^x = (m^x_{ij})_{i,j\in \I}$ with
\begin{align}\label{eq:coxeter-datum}
s((\sigma^x_i\sigma_j)^{m^x_{ij}}) &= x,& i,j&\in \I,& x&\in \cX.
\end{align}
The \emph{Coxeter groupoid} $\cW(\cX, \rho, \bM)$ 
associated to $(\cX, \rho, \bM)$ \cite[Definition 1]{HY} 
is the groupoid  presented by generators $\Q_{\rho}$ with relations
\begin{align}\label{eq:def-coxeter-gpd}
 (\sigma_i^x\sigma_j)^{m^x_{ij}} &= \id_x, &i, j\in \I,\,  &x\in \cX.
\end{align}

\smallbreak\subsubsection{Generalized root system}\label{subsubsec:grs} 
A generalized root system (GRS for short) is a collection $\cR:= (\cX, \rho, \cC, \Delta)$, where  
$\cC = (C^x)_{x\in \cX}$ is a family of generalized Cartan matrices $C^x = (c^x_{ij})_{i,j \in \I}$, cf.  \cite{K-libro},
and  $\Delta = (\Delta^x)_{x\in \cX}$ is a family of subsets $\Delta^x \subset \mathbb Z^{\I}$. We need the following notation: 
Let $\{\alpha_i\}_{i\in\I}$ be the  canonical basis of $\Z^{\I}$ and
define $s_i^x\in GL(\Z^{\I})$ by
$s_i^x(\alpha_j) = \alpha_j-c_{ij}^x\alpha_i$,  $i, j \in \I$, $x \in \cX$.
The collection should satisfy the following axioms:
\begin{align}
c^x_{ij}&=c^{\rho_i(x)}_{ij} &  \mbox{for all }&x \in \cX, \, i,j \in \I. \\
\label{eq:def root system 1}
\Delta^x &= \Delta^x_+ \cup \Delta^x_-, & \Delta^x_{\pm} &:= \pm(\Delta^x \cap \mathbb N_0^{\I}) \subset \pm\mathbb N_0^{\I};
\\ \label{eq:def root system 2}
\Delta^x \cap \mathbb Z \alpha_i &= \{\pm \alpha_i \};& &
\\\label{eq:def root system 3}
s_i^x(\Delta^x)&=\Delta^{\rho_i(x)};& 
\\ \label{eq:def root system 4}
(\rho_i\rho_j)^{m_{ij}^x}(x)&=(x), & m_{ij}^x &:=|\Delta^x \cap (\mathbb N_0\alpha_i+\mathbb N_0 \alpha_j)|,
\end{align}
for all $x \in \cX$, $i \neq j \in \I$.
We call $\Delta^x_+$, respectively $ \Delta^x_-$, the set of \emph{positive}, respectively \emph{negative}, roots.
Let $\cG = \cX \times GL_{\theta}(\Z) \times \cX$, $\varsigma_i^x = (x, s_i^x,\rho_i(x))$, $i \in \I$, $x \in \cX$,
and $\cW = \cW(\cX, \rho, \cC)$ the subgroupoid of $\cG$ generated by all the $\varsigma_i^x$,
i.e., by the image of the morphism of quivers $\Q_{\rho} \to \cG$, $\sigma_i^x \mapsto \varsigma_i^x$. 
There is a Coxeter matrix $\bm^x = (m^x_{ij})_{i,j\in \I}$, where $m^x_{ij}$ is the smallest natural number such that
$(\varsigma_i^x\varsigma_j)^{m^x_{ij}} = \id x$. Then  
 $\bM = (\bm^x)_{x\in \cX}$ fits into a Coxeter datum $(\cX, \rho, \bM)$, and
there is an isomorphism  of groupoids 
$\xymatrix{ \cW(\cX, \rho, \bM) \ar@{->>}[r] &\cW = \cW(\cX, \rho, \cC)}
$ \cite{HY}; this is called the \emph{Weyl groupoid} of $\cR$.
If $w \in \cW(x, y)$, then $w(\Delta^x)= \Delta^y$, by \eqref{eq:def root system 3}.
The sets of \emph{real} roots at $x \in \cX$ are
$(\Delta^{\re})^x = \bigcup_{y\in \cX}\{ w(\alpha_i): \ i \in \I, \ w \in \cW(y,x) \}$; correspondingly the \emph{imaginary} roots are 
$(\Delta^{\im})^x = \Delta^{x} - (\Delta^{\re})^x$.
Assume that $\cW$ is connected.
Then the following conditions are equivalent
\cite[Lemma 2.11]{CH-at most 3}:
\begin{itemize} \itemsep1pt \parskip0pt \parsep0pt
	\item $\Delta^x$ is finite for some $x\in \cX$,
	\item $\Delta^x$ is finite for all $x\in \cX$,
	\item $(\Delta^{\re})^x$ is finite for all $x\in \cX$,
\item $\cW$ is finite. 
\end{itemize}
If these hold, then all roots are real \cite{CH-at most 3};  we say that $\cR$ is \emph{finite}.
We now discuss two examples of GRS, central for the subsequent discussion.

\begin{example}\label{exa:grs-lie-super} \cite{AA}
Let $\ku$ be a field of characteristic $\ell \geq 0$, $\theta\in\N$, $\bp\in\G_2^\theta$ and $A =(a_{ij})\in\ku^{\theta\times\theta}$.
We assume $\ell \neq 2$ for simplicity. Let $\h = \ku^{2\theta-\rk A}$.  
Let $\g(A,\bp)$ be the Kac-Moody Lie superalgebra over $\ku$ defined as in \cite{K-libro}; 
it is generated  by $\h$, $e_i$ and $f_i$, $i\in \I$, and the parity is given by
$|e_i|=|f_i|=p_i$,  $i\in\I$, $|h|=0$, $h\in\h$. Let $\Delta^{A, \bp}$ be the root system of $\g(A,\bp)$.
We make the following technical assumptions:
\begin{align}\label{eq:symmetrizable}
a_{jk} &=0 \implies a_{kj}=0, & j&\neq k; \\
\label{eq:f-nilpotent}
\ad f_i &\text{ is locally nilpotent in } \g(A,\bp), & i&\in\I.
\end{align}
The matrix $A$ is \emph{admissible} if \eqref{eq:f-nilpotent} holds \cite{Ser-superKM}.
Let $C^{A,\bp}= (c_{ij}^{A,\bp})_{i,j\in\I}$ be given by 
\begin{align}\label{eq:cryst-g(A)}
c_{ij}^{A,\bp}&:=-\min\{m\in\N_0:(\ad f_i)^{m+1} f_j= 0 \}, i\neq j \in \I, & c_{ii}^{A,\bp}&:= 2.
\end{align}
We need the following elements of $\ku$:

\begin{align}
\label{eq:def dn, i par} \mbox{if } p_i &= 0, & d_m&= m\, a_{ij}+\binom{m}{2} a_{ii};  \\
\label{eq:def dn, i impar} \mbox{if }p_i &= 1; &
d_m&= \left\{ \begin{array}{ll} k\, a_{ii}, & m=2k, \\ k\, a_{ii}+a_{ij},& m=2k+1; \end{array}\right.  
\\
\label{eq:formula nu}
\nu_{j,0}&=1, & \nu_{j,n}&= \prod_{t=1}^n (-1)^{p_i((t-1)p_i+p_j)}d_t;        \\
\label{eq:formula mu} 
\mu_{j,0}&=0, & \mu_{j,n}&= (-1)^{p_ip_j}n\left(\prod_{t=2}^n (-1)^{p_i((t-1)p_i+p_j)}d_t\right)a_{ji} .       
\end{align}

With the help of these scalars, we define a reflection $r_i(A,\bp) = (r_iA, r_i\bp)$, where
$r_i\bp=(\overline{p}_j)_{j\in\I}$, with $\overline{p}_j=p_j-c_{ij}^{A,\bp}p_i$, and $r_i A=(\overline{a}_{jk})_{j,k\in\I}$, with
\begin{align}\label{eq:formula matriz si(A)}
\overline{a}_{jk}=\left\{\begin{array}{ll} -c_{ik}^{A,\bp}\mu_{j,-c_{ij}^{A,\bp}}a_{ii}+\mu_{j,-c_{ij}^{A,\bp}}a_{ik} & \\
\qquad -c_{ik}^{A,\bp}\nu_{j,-c_{ij}^{A,\bp}}a_{ji}+\nu_{j,-c_{ij}^{A,\bp}}a_{jk}, & j,k\neq i; \\
c_{ik}^{A,\bp}a_{ii}-a_{ik}, & j=i\neq k; \\ -\mu_{j,-c_{ij}^{A,\bp}}a_{ii}-\nu_{j,-c_{ij}^{A,\bp}}a_{ji}, & j\neq k=i; \\ a_{ii}, & j=k=i.\end{array} \right.
\end{align}

\begin{theorem}\label{thm:isomorfismo Ti} There is an isomorphism $T_i^{A,\bp}:\g(r_i (A,\bp))\to\g(A,\bp)$
of Lie superalgebras  given (for an approppriate basis $(h_i)$ of $\h$) by
\end{theorem}
\begin{equation}\label{eq:def Ti}
\begin{split}
T_i^{A,\bp}(e_j)&= \begin{cases} 
(\ad e_i)^{-c_{ij}^{A,\bp}} (e_j), & i\neq j\in\I, 
\\ f_i, & j=i 
\end{cases}  
\\
T_i^{A,\bp}(f_j)&= \begin{cases} (\ad f_i)^{-c_{ij}^{A,\bp}}f_j, & j\in\I,j\neq i, \\
(-1)^{p_i}e_i, & j=i,\end{cases}   
\\
T_i^{A,\bp}(h_j)&= \begin{cases} \mu_{j,-c_{ij}^{A,\bp}} h_i+\nu_{j,-c_{ij}^{A,\bp}} h_j,& i\neq j\in\I \\
-h_i, & j=i,\\ h_j, & \theta + 1 \leq j \leq 2\theta-\rk A.\end{cases} 
\end{split}
\end{equation}
Assume that $\dim \g(A,\bp) < \infty$; then \eqref{eq:symmetrizable} and \eqref{eq:f-nilpotent} hold. Let 
\begin{align*}
 \cX =&\{ r _{i_1}\cdots r _{i_n}(A,\bp)\,|\, n\in \N _0, i_1,\ldots ,i_n\in \I\}.
 \end{align*}
Then $(\cX, r, \cC,  \Delta)$,  where $\cC = (C^{(B, \bq)})_{(B, \bq)\in \cX}$ and $\Delta = (\Delta^{(B, \bq)})_{(B, \bq)\in \cX}$,  
is a finite GRS, an invariant of 
$\g(A,\bp)$. \end{example}

\begin{example}\label{exa:grs-nichols-ss-yd} Let $H$ be a Hopf algebra, assumed semisimple for easiness. 
Let $M\in \ydh$ be finite-dimensional,
with a fixed decomposition $M =  M_1\oplus \dots \oplus M_{\theta}$, where 
$M_1,\dots,M_{\theta} \in \Irr \ydh$. Then $T(M)$ and $\toba(M)$ are $\zt$-graded, by $\deg x =\alpha_i$ 
for all $x\in M_i$, $i\in \I_\theta$. Recall that $\zt_{\ge 0} = \sum_{i\in \I_\theta}\Z_{\ge 0}\alpha_i$.

\begin{theorem}\label{th:grs-semisimple} \emph{\cite{HeckSch, HS-israel}} 
If $\dim \toba(M) < \infty$, then  $M$ has a finite GRS.
\end{theorem}

We discuss the main ideas of the proof. Let $i\in \I = \I_\theta$. We define
$M'_i=V_i^*$, 
\begin{align*}
c_{ij}^M &= - \sup \{h\in \N_0: \ad_c^h (M_i)(M_j)\neq 0 \text{ in }\toba(M)\},& i&\neq j, \quad  c^M_{ii} = 2;
\\  M'_j&= \ad_c^{-c_{ij}} (M_i)(M_j), & \rho_i(M) &= M'_1 \oplus \dots \oplus M'_\theta.
\end{align*}
 Then $\dim \toba(M) = \dim \toba(\rho_i(M))$ and  $C^M = (c^M_{ij})_{i,j\in\I}$ is a
generalized Cartan matrix  \cite{AHS}. 
Also, $M'_j$ is irreducible \cite[3.8]{AHS},  \cite[7.2]{HeckSch}.
Let  $\cX$ be the set of  objects in $\ydh$ with fixed decomposition  (up to isomorphism)  of the form
\begin{align*}
 \{ \rho _{i_1}\cdots \rho _{i_n}(M)\,|\, n\in \N _0, i_1,\ldots ,i_n\in \I\}.
 \end{align*}
Then $(\cX, \rho, \cC)$,  where $\cC = (C^N)_{N\in \cX}$,  is a crystallographic datum. Next we need:

\smallbreak\noindent $\bullet$ \cite[Theorem 4.5]{HeckSch}; \cite{HeckGran07} There exists a totally ordered index set $(L,\le)$ and families 
$(W_l)_{l\in L}$ in $\Irr\ydh$, $(\beta_l)_{l\in L}$   such that $\toba (M)\simeq \otimes _{l\in L}\toba (W_l)$ as 
$\zt$-graded objects in $\ydh$, where $\deg x =\beta_l$ for all $x\in W_l$, $l\in L$.  

\smallbreak
Let $\Delta_{\pm}^M = \{\pm\beta_l: l\in L\}$, $\Delta^M = \Delta _+^M \cup \Delta_-^M$, $\Delta = (\Delta^N)_{N\in \cX(M)}$.
Then $\cR = (\cX, \rho, \cC, \Delta)$ is a finite GRS.
\end{example}

\begin{theorem}\label{th:classif-GRS}  \emph{\cite{CH}}
The classification of all  finite GRS is known.
\end{theorem}

The proof is a combinatorial tour-de-force and requires computer calculations.
It is possible to recover from this result the classification of the finite-dimensional contragredient Lie superalgebras 
in arbitrary characteristic \cite{AA}. However, the list of \cite{CH} is substantially larger than the classifications of 
the alluded Lie superalgebras or the braidings of diagonal type with finite-dimensional Nichols algebra.

\smallbreak\subsection{Nichols algebras of diagonal type}\label{subsec:nichols-diagonal}
Let $G$ be a finite group. We denote $\ydgg = \ydh$ for $H = \kk G$. So
$M\in \ydgg$ is a left $G$-module with a $G$-grading  $M = \oplus_{g\in G}M_g$ such that
$t\cdot M_g = M_{tgt^{-1}}$, for all $g,t\in G$.
If $M,N \in \ydgg$,  then the braiding $c: M\otimes N \to N\otimes 
M$ is given by $c(m\otimes n) = g\cdot n \otimes m$, $m\in M_g$, $n\in N$, $g\in G$.
Now assume that $G = \Gamma$ is a finite abelian group. Then every $M\in \ydg$ is a $\Gamma$-graded $\Gamma$-module, hence of the form
$M =   \oplus _{g\in \Gamma, \chi\in \VGamma} M_g^{\chi}$, where $M_g^{\chi}$ is the $\chi$-isotypic component of $M_g$. 
So $\ydg$ is just the category
of $\Gamma \times \VGamma$-graded modules, with the braiding  $c: M\otimes N \to N\otimes 
M$  given by $c(m\otimes n) = \chi(g) n \otimes m$, $m\in M^\eta_g$, $n\in N_t^\chi$, $g, t\in G$, $\chi, \eta\in \VGamma$.
Let $\theta\in \N$, $\I = \I_{\theta}$.

\begin{definition}\label{def:diagonal} Let $\bq = (q_{ij})_{i, j\in \I}$ be a matrix with entries in $\kk^{\times}$.
A braided vector space $(V,c)$ is of \emph{diagonal type} with matrix $\bq$ if $V$ has a basis $(x_i)_{i\in \I}$
with
\begin{align}\label{eq:diagonal}
c(x_i \otimes x_j) &= q_{ij} x_j \otimes x_i, & i, j &\in \I.
\end{align}
\end{definition}
Thus, every finite-dimensional $V\in \ydg$ is a braided vector space of diagonal type. 
Question \ref{item:method-nichols}, more precisely (a$_1$), has a complete answer in this setting.
First we can assume that $q_{ii} \neq 1$ for $i\in \I$, as otherwise $\dim \toba(V) = \infty$. 
Also, let $\bq' = (q'_{ij})_{i, j\in \I} \in (\kk^{\times})^{\I \times \I}$ and $V'$ a braided vector space with matrix $\bq'$.
If $q_{ii} = q'_{ii}$ and $q_{ij}q_{ji} = q'_{ij}q'_{ji}$ for all  $j\neq i\in \I_\theta$, then $\toba(V)\simeq \toba(V')$
as braided vector spaces.

\begin{theorem}\label{th:classif-diagonal}  \emph{\cite{H-classif}}
The classification of all braided vector spaces of diagonal type with finite-dimensional Nichols algebra is known.
\end{theorem}

The proof relies on the Weyl groupoid introduced in \cite{H-inv}, a particular case of Theorem \ref{th:grs-semisimple}.
Another fundamental ingredient is the following result, generalized at various levels in  \cite{HeckGran07, HeckSch, HS-israel}.

\begin{theorem}\label{th:pbw}  \emph{\cite{Khar99}}
Let $V$ be a braided vector space of diagonal type. Every Hopf algebra quotient of $T(V)$ has a PBW basis.
\end{theorem}

The classification in Theorem \ref{th:classif-diagonal} can be organized as follows: 

\smallbreak \noindent $\diamond$ For most of the matrices $\bq = (q_{ij})_{i, j\in \I_\theta}$ in the list of \cite{H-classif} 
there is a field $\ku$ and a pair $(A, \bp)$ as in Example \ref{exa:grs-lie-super} such that $\dim \g(A, \bp) < \infty$, and $\g(A, \bp)$ 
has the same GRS as the Nichols algebra corresponding to $\bq$ \cite{AA}.

\smallbreak \noindent $\diamond$ Besides these, there are 12 (yet) unidentified examples.

\noindent We believe that Theorem \ref{th:classif-diagonal} can be proved from Theorem \ref{th:classif-GRS}, via Example \ref{exa:grs-lie-super}.

\begin{theorem}\label{th:presentation-diagonal}  \emph{\cite{A-jems, A-crelle}}
An efficient set of defining relations of each finite-dimen\-sional Nichols algebra of a braided vector space of diagonal type is known.
\end{theorem}

The proof uses most technical tools available in the theory of Nichols algebras; of interest in its own is the introduction of the notion of convex order 
in Weyl groupoids. As for other classifications above, it is not possible to state precisely the list of relations. 
We just mention different types of relations that appear.

\begin{itemize}\renewcommand{\labelitemi}{$\circ$}\itemsep1pt \parskip0pt \parsep0pt
\item Quantum Serre relations, i.e., $\ad_c(x_i)^{1- a_{ij}} (x_j)$ for suitable $i\neq j$.
\item Powers of root vectors, i.e., $x_\beta^{N_\beta}$, where the $x_\beta$'s are part of the PBW basis.
\item More exotic relations; they involve 2, 3, or at most 4 $i$'s in $\I$.
\end{itemize}

\smallbreak\subsection{Nichols algebras of rack type}\label{subsec:nichols-rack}
We now consider Nichols algebras of objects in $\ydgg$, where $G$ is a finite not necessarily abelian group. 
The category $\ydgg$ is semisimple and the simple objects are parametrized by pairs $(\Oc, \rho)$, where $\Oc$ is a conjugacy class
in $G$ and $\rho\in \Irr G^x$, for a fixed $x\in \Oc$; the corresponding simple Yetter-Drinfeld module 
$M(\Oc, \rho)$ is $\Ind_{G^x}^G \rho$ as a module. The braiding $c$ is described in terms of the conjugation in $\Oc$. To describe the related suitable class, 
we recall that a rack is a set $X \neq \emptyset$  with a map $\trid: X \times X \to X$ satisfying 
\begin{itemize}\renewcommand{\labelitemi}{$\circ$}\itemsep1pt \parskip0pt \parsep0pt
 \item  $\varphi_x := x\trid \underline{\quad}$ is a
bijection for every $x \in X$.  
\item $x\trid (y\trid z) = (x\trid y) \trid (x\trid z)$ for all $x,y,z \in X$ (self-distributivity).
\end{itemize}

For instance, a conjugacy class $\Oc$ in $G$ with the operation $x\trid y = xyx^{-1}$, $x, y \in \Oc$
is a rack; actually we only consider  racks realizable as conjugacy classes.
Let $X$ be a rack and  $\mathfrak X = (X_k)_{k\in I}$ a decomposition of $X$, i.e., a disjoint family of subracks with $X_l\trid X_k = X_k$ for 
all $k,l\in I$. 

\begin{definition}\label{def:braiding-rack}\cite{AG-adv} A \emph{2-cocycle of degree $\n = (n_k)_{k\in I}$},
associated to $\mathfrak X$, is a family $\bq
= (q_k)_{k\in I}$ of maps $q_k: X \times X_k \to \GL(n_k,\kk)$
such that
\begin{align}\label{eqn:non-ppal-cocycle-braiding}
 q_k(i,j\trid h)q_k(j,h)&= q_k(i\trid j,i\trid h)q_k(i,h),& i,j &\in X, \, h\in X_k, \, k\in I.
\end{align}
Given such  $\bq$, let $V = \oplus_{k\in I} \kk X_k\otimes\kk^{n_k}$ and let $c^{\bq}\in \GL(V\otimes V)$ be given by
\begin{align*}
c^{\bq}(x_i v\otimes x_jw) &= x_{i\trid j}q_k(i,j)(w)\otimes x_iv,& i&\in X_l, \, j\in X_k, \, v\in\kk^{n_l}, \, w\in\kk^{n_k}.
\end{align*}
Then $(V, c^{\bq})$ is a braided vector space called of \emph{rack type}; its Nichols algebra is denoted $\toba(X, \bq)$.
If $\mathfrak X = (X)$, then we say that $\bq$ is \emph{principal}.
\end{definition}

Every finite-dimensional $V\in \ydgg$ is a braided vector space of rack type \cite[Theorem 4.14]{AG-adv}. 
Question  (a$_1$) in this setting has partial answers in three different lines: computation of some finite-dimensional Nichols
algebras, Nichols algebras of reducible Yetter-Drinfeld modules and collapsing of racks.

\smallbreak\subsubsection{Finite-dimensional Nichols algebras of rack type}\label{subsubsec:nichols-rack-findim}
The algorithm to compute a Nichols algebra $\toba(V)$  is as follows: 
compute the space $\J^i(V) = \ker \Q_i$ of relations of degree $i$, for $i = 2, 3, \dots, m$; then compute 
the $m$-th partial Nichols algebra $\widehat{\toba_m}(V) = T(V)/ \langle \oplus_{2\le i \le m} \J^i(V)\rangle$,
say with a computer program. If lucky enough to get $\dim \widehat{\toba_m}(V) < \infty$, then check whether it is a Nichols algebra, e.g. via skew-derivations;
otherwise go to $m+1$. The description of $\J^2(V) = \ker (\id + c)$ is not difficult \cite{GG} but for higher degrees it turns out to be very complicated.
We list all known examples of finite-dimensional Nichols algebras $\toba(X, \bq)$ with $X$ indecomposable and $\bq$ principal 
and abelian ($n_1 = 1$).

\begin{example}\label{exa:nichols-sim} Let $\Oc^m_d$ be the conjugacy class of $d$-cycles in
$\sm$, $m \ge 3$. We start with the rack of transpositions in $\sm$ and the cocycles $-1$, $\chi$ that arise from the $\rho \in \Irr \sm^{(12)}$
with $\rho(12) = -1$, see \cite[(5.5), (5.9)]{MS}. 
Let $V$ be a vector space with basis $(x_{ij})_{(ij)\in \Oc^m_2}$ and consider the relations 
\begin{align}
 \label{eq:trasp-quadratic}
 x_{ij}^2 &= 0, & &(ij)\in \Oc^m_2;
 \\   \label{eq:trasp-far}
 x_{ij}x_{kl} + x_{kl}x_{ij} &= 0, & &(ij), (kl)\in \Oc^m_2, \,  \vert \{i,j,k,l\}\vert = 4;
 \\  \label{eq:trasp-far-sign}
 x_{ij}x_{kl} - x_{kl}x_{ij} &= 0, & &(ij), (kl)\in \Oc^m_2, \,  \vert \{i,j,k,l\}\vert = 4;
 \\ \label{eq:trasp-braided}
 x_{ij}x_{ik} + x_{jk}x_{ij}+ x_{ik}x_{jk} &= 0, & &(ij), (ik), (jk)\in \Oc^m_2, \, \vert \{i,j,k\}\vert = 3;
 \\ \label{eq:trasp-braided-sgn}
x_{ij}x_{ik} - x_{jk}x_{ij} - x_{ik}x_{jk} &= 0, & &(ij), (ik), (jk)\in \Oc^m_2, \, \vert \{i,j,k\}\vert = 3.
\end{align}
The quadratic algebras $\toba_m := \widehat{\toba_2}( \Oc^m_2, -1) = T(V)/\langle \text{\eqref{eq:trasp-quadratic}, \eqref{eq:trasp-far}, \eqref{eq:trasp-braided}}\rangle$
and
$\Ee_m:= \widehat{\toba_2}( \Oc^m_2, \chi) = T(V)/\langle \text{\eqref{eq:trasp-quadratic}, \eqref{eq:trasp-far-sign}, \eqref{eq:trasp-braided-sgn}}\rangle$
were considered in \cite{MS}, \cite{FK} respectively; $\Ee_m$ are named the \emph{Fomin-Kirillov algebras}. It is known that

\smallbreak \noindent $\circ$ The Nichols algebras $\toba( \Oc^m_2, -1)$ and $\toba( \Oc^m_2, \chi)$ are twist-equivalent, hence have the same Hilbert series. 
Ditto for the  algebras $\toba_m$ and $\Ee_m$  \cite{V}.
 
\smallbreak \noindent $\circ$ If $3 \le m \le 5$, then $\toba_m = \toba( \Oc^m_2, -1)$ and $\Ee_m =\toba( \Oc^m_2, \chi)$ are finite-dimensional  
\cite{FK, MS, GG} (for $m=5$ part of this was done by Gra\~na). In fact
\begin{align*}
\dim \toba_3 &= 12, &\dim \toba_4 &= 576, &\dim \toba_5 &= 8294400.
\end{align*}
But for $m\geq 6$, it is not known whether the Nichols algebras $\toba( \Oc^m_2, -1)$ and $\toba( \Oc^m_2, \chi)$ have finite
dimension or are quadratic.
\end{example}

\begin{example}\label{exa:nichols-cube} \cite{AG-adv}
The Nichols algebra $\toba( \Oc^4_4, -1)$ is quadratic, has the same Hilbert series as $\toba( \Oc^4_2, -1)$ 
and is generated by $(x_\sigma)_{\sigma\in \Oc_4^4}$ with defining relations
\begin{align}
 x_\sigma^2 &= 0, & &\\ 
 x_\sigma x_{\sigma^{-1}}+x_{\sigma^{-1}} x_\sigma &= 0, & & \\ 
x_\sigma x_\kappa+x_\nu x_\sigma+x_\kappa x_\nu&= 0,&
\sigma\kappa &=\nu\sigma, \ \kappa\neq \sigma\neq \nu\in \Oc_4^4.
\end{align}
\end{example}

\begin{example}\label{exa:nichols-affines} \cite{grania1}  Let $A$ be a finite abelian group
and $g\in \Aut A$. The  \emph{affine} rack $(A, g)$ is the set $A$ with product  $a\trid b = g(b) + (\id - g) (a)$,  $a, b \in A$. 
Let $p\in \N$ be a prime, $q =p^{v(q)}$ a power of $p$,  $A = \F_q$ and $g$ the  multiplication by  $N \in \F_q^{\times}$; 
let $X_{q, N} = (A, g)$. Assume that $q = 3$, 4, 5, or 7, with
$N = 2$, $\omega\in \F_4 - \F_2$, 2 or 3, respectively. Then 
$\dim \toba(X_{q, N}, -1) = q \varphi(q)(q-1)^{q-2}$, $\varphi$ being the Euler function, and
$\J(X_{q, N}, -1) = \langle\J^2 + \J^{v(q)(q-1)}\rangle$, where  $\J^2$ is generated by 
\begin{align}
&x_i^2, &&\text{always}\\
&x_ix_j+x_{-i+2j}x_i+x_jx_{-i+2j}, &&\text{for $q=3$,}\\
&x_ix_j+x_{(\omega+1)i+\omega j}x_i+x_jx_{(\omega+1)i+\omega j,} &&\text{for $q=4$,}\\
&x_ix_j+x_{-i+2j}x_i+x_{3i-2j}x_{-i+2j}+x_jx_{3i-2j}, &&\text{for $q=5$,}\\
&x_ix_j+x_{-2i+3j}x_i+x_jx_{-2i+3j}, &&\text{for $q=7$,}
\end{align}
with $i, j \in \F_q$; and $\J^{v(q)(q-1)}$ is generated by $\displaystyle\sum_{h}T^h(V)\J^2T^{v(q)(q-1) -h -2}(V)$ and 
\begin{align}
&(x_\omega x_1x_0)^2 + (x_1x_0x_\omega)^2  + (x_0x_\omega x_1)^2, &&\text{for $q=4$,}  \\
&(x_1x_0)^2 + (x_0x_1)^2, &&\text{for $q=5$,} \\
&(x_2x_1x_0)^2 + (x_1x_0x_2)^2 + (x_0x_2x_1)^2, &&\text{for $q=7$.}
\end{align}
Of course $X_{3, 2} = \Oc_2^3$; also $\dim \toba(X_{4, \omega}, -1) = 72$. By duality, we get
\begin{align*}
\dim \toba(X_{5, 3}, -1) = \dim \toba(X_{5, 2}, -1) &= 1280, \\\dim \toba(X_{7, 5}, -1)  = \dim \toba(X_{7, 3}, -1) &= 326592.
\end{align*}
\end{example}

\begin{example}\label{exa:nichols-tetra}\cite{HLV} There is another finite-dimensional Nichols algebra associated to $X_{4, \omega}$
with a cocycle $\bq$ with values $\pm \xi$, where $1\neq \xi\in\G_3$. Concretely, $\dim \toba(X_{4, \omega}, \bq) = 5184$ and
$\toba(X_{4, \omega}, \bq)$ can be presented by generators $(x_i)_{i\in \F_4}$ with defining relations
\begin{gather*}
\begin{aligned}
x_{0}^{3}=x_{1}^{3}=x_{\omega}^{3}=x_{\omega^2}^{3}&=0, && \\
\xi ^2x_{0}x_{1} + \xi  x_{1}x_{\omega} - x_{\omega}x_{0} &=0,&  \xi ^2x_{0}x_{\omega} + \xi  x_{\omega}x_{\omega^2} - x_{\omega^2}x_{0} &= 0,\\
\xi  x_{0}x_{\omega^2} -\xi ^2x_{1}x_{0} + x_{\omega^2}x_{1} &=0, & \xi  x_{1}x_{\omega^2} + \xi ^2x_{\omega}x_{1} + x_{\omega^2}x_{\omega} &=0,
\end{aligned}
\\
x_{0}^2x_{1}x_{\omega}x_{1}^2 + x_{0}x_{1}x_{\omega}x_{1}^2x_{0} + x_{1}x_{\omega}x_{1}^2x_{0}^2 + x_{\omega}x_{1}^2x_{0}^2x_{1} + x_{1}^2x_{0}^2x_{1}x_{\omega} + x_{1}x_{0}^2x_{1}x_{\omega}x_{1}  \\
+ x_{1}x_{\omega}x_{1}x_{0}^2x_{\omega} + x_{\omega}x_{1}x_{0}x_{1}x_{0}x_{\omega} + x_{\omega}x_{1}^2x_{0}x_{\omega}x_{0} =0.
\end{gather*}

\end{example}

\smallbreak\subsubsection{Nichols algebras of decomposable Yetter-Drinfeld modules over groups}\label{subsubsec:nichols-rack-decomposable}
The ideas of Example \ref{exa:grs-nichols-ss-yd} in the context of decomposable Yetter-Drinfeld modules over groups were pushed further 
in a series of papers culminating with a remarkable classification result \cite{HV2}. Consider the groups
\begin{align}
\Gamma_n &=  \langle a,b,\nu \vert ba=\nu ab,\quad \nu a=a\nu^{-1},\quad \nu b=b\nu,\quad \nu ^n=1\rangle,\quad n\ge 2;
\\ T &= \langle\zeta ,\chi _1,\chi _2\vert   \zeta \chi_1=\chi_1\zeta ,\quad
  \zeta \chi_2=\chi_2\zeta ,\quad   \chi_1\chi_2\chi_1=\chi_2\chi_1\chi_2,\quad
  \chi_1^3=\chi_2^3\rangle.
\end{align}

\smallbreak \noindent $\circ$ \cite{HS2} Let $G$ be a quotient of $\Gamma_{2}$. Then there exist $V_1$, $W_1 \in \Irr\ydgg$ such that 
$\dim V_1 = \dim W_1 = 2$ and $\dim\toba(V_1\oplus W_1) = 64 = 2^6$.

\smallbreak \noindent $\circ$ \cite{HV2} Let $G$ be a quotient of $\Gamma_{3}$. Then there exist $V_2$, $V_3$, $V_4$, $W_2, W_3, W_4 \in \Irr\ydgg$ such that 
$\dim V_2 = 1$, $\dim V_3 = \dim V_4 = 2$, $\dim W_2 = \dim W_3 = \dim W_4 = 3$ and $\dim\toba(V_2\oplus W_2) = \dim\toba(V_3\oplus W_3) = 10368  =  2^73^4$,  
$\dim\toba(V_4\oplus W_4) = 2304 = 2^{18}$.

\smallbreak \noindent $\circ$ \cite{HV1} Let $G$ be a quotient of $\Gamma_{4}$. Then there exist $V_5$, $W_5 \in \Irr \ydgg$ such that 
$\dim V_5 = 2$, $\dim W_5 = 4$ and $\dim\toba(V_5\oplus W_5) =  262144 = 2^{18}$.

\smallbreak \noindent $\circ$ \cite{HV1} Let $G$ be a quotient of $T$. Then there exist $V_6$, $W_6 \in \Irr\ydgg$ such that 
$\dim V_6 = 1$, $\dim W_6 = 4$ and $\dim\toba(V_6\oplus W_6) = 80621568 = 2^{12}3^9$.

\begin{theorem} \emph{\cite{HV2}} 
Let $G$ be a non-abelian group and $V, W\in \Irr \ydgg$   such that $G$ is generated by the support of $V\oplus W$. 
Assume that     $c^2_{\vert V\otimes W} \neq \id$  and that $\dim \toba (V\oplus W) < \infty$. 
Then $V\oplus W$ is one of $V_i\oplus W_i$, $i\in \I_6$, above, and correspondingly $G$ is a quotient of
either $\Gamma_n$, $2\le n\le 4$, or $T$.
\end{theorem}

\smallbreak\subsubsection{Collapsing racks}\label{subsubsec:nichols-rack-collapsing}
Implicit in Question  (a$_1$) in the setting of racks is the need  to compute all non-principal 2-cocycles for a fixed rack $X$.
Notably, there exist criteria that dispense of this computation. To state them and explain their significance, we need some terminology.
All racks below are finite.

\smallbreak \noindent $\circ$ A rack $X$ is \emph{abelian} when $x\trid y = y$, for all $x, y \in X$.

\smallbreak \noindent $\circ$
A rack is \emph{indecomposable} when it is not a disjoint union
of two proper subracks.

\smallbreak \noindent $\circ$ A rack $X$ with $\vert X\vert > 1$ is \emph{simple} when for any projection of racks
$\pi: X\to Y$, either $\pi$ is an isomorphism or $Y$ has only one element. 

\begin{theorem} \emph{\cite[3.9, 3.12]{AG-adv}, \cite{jo}}
 Every simple rack is isomorphic to one of:

\begin{enumerate}\renewcommand{\theenumi}{(\alph{enumi})}\itemsep1pt \parskip0pt \parsep0pt
    \item Affine racks $(\F_p^t, T)$, where $p$ is a prime, $t\in \N$,  and
    $T$ is the companion matrix of a monic irreducible polynomial $f\in \F_p[\mathtt X]$ of degree $t$,
     $f\neq \mathtt X, \mathtt X-1$.

    \item\label{item:twisted-conjugacy} Non-trivial (twisted) conjugacy classes in simple groups.

    \item\label{item:twisted-homogeneous} Twisted conjugacy classes of type $(G,u)$, where $G = L^t$, with $L$ a simple non-abelian group and $1 < t\in \N$;
and $u\in \Aut (L^t)$ acts by $u(\ell_1,\ell_2, \dots,\ell_t)=(\theta(\ell_t),\ell_1,\ell_2, \dots,\ell_{t-1})$, 
where $\theta\in\Aut(L)$.   
\end{enumerate}
\end{theorem}

\begin{definition}\label{def:rack-typeD}
 \cite[3.5]{AFGV-ampa} We say that a finite rack $X$ is \emph{of type D}
when there are a decomposable subrack
$Y = R\coprod S$, $r\in R$ and  $s\in S$ such that
$r\trid(s\trid(r\trid s)) \neq s$. 

Also,
$X$ is \emph{of type F} \cite{ACG-I} if there are a disjoint family of subracks $(R_a)_{a \in \I_4}$
and a family $(r_a)_{a \in \I_4}$ with $r_a\in R_a$, such that 
$R_a \triangleright R_b = R_b$,
$r_a\triangleright r_b \neq r_b$, for all $a\neq b \in \I_4$.
\end{definition}

An indecomposable rack $X$ \emph{collapses} when $\dim \toba(X, \bq) = \infty$
for every finite \emph{faithful} 2-cocycle $\bq$ (see  \cite{AFGV-ampa} for the definition of faithful).

\begin{theorem}\label{th:type4}  \emph{\cite[3.6]{AFGV-ampa}; \cite[2.8]{ACG-I}} If a rack is of type D or F, then it collapses.
\end{theorem}

The proofs use results on Nichols algebras from \cite{AHS, HeckSch, CH}.

If a rack projects onto a rack of type D (or F), then it is also
of type D (or  F), hence it collapses by Theorem \ref{th:type4}.
Since every indecomposable rack $X$, $\vert X\vert > 1$, projects onto a simple rack,
it is natural to ask for the determination of all \emph{simple} racks of type D or F.
A rack is \emph{cthulhu}
if it is neither  of type D nor F;
it is \emph{sober} if  every subrack is either abelian or indecomposable \cite{ACG-I}. Sober implies cthulhu.

\smallbreak \noindent $\circ$  Let $m\geq 5$. Let $\Oc$ be either $\Oc^{\sm}_{\sigma}$, if 
$\sigma\in \sm -\am$, or else $\Oc^{\am}_{\sigma}$ if $\sigma\in \am$. The type of $\sigma$ is formed by the lengths of the cycles in its decomposition.

\smallbreak \noindent $\diamond$ \cite[4.2]{AFGV-ampa} If the type of $\sigma$ is $(3^2)$, $(2^2,3)$, $(1^n,3)$, $(2^4)$, $(1^2,2^2)$, $(2,3)$,  $(2^3)$, or 
$(1^n,2)$, then $\Oc$ is cthulhu.
 If the type of $\sigma$ is $(1,2^2)$, then $\Oc$ is sober. 
 
\smallbreak \noindent $\diamond$ \cite{F} Let $p\in \N$ be a prime. Assume the type of $\sigma$ is $(p)$. 
 If $p = 5, 7$ or  not of the form $(r^k - 1)/(r -1)$,  $r$  a prime power,  then $\Oc$ is sober; otherwise $\Oc$ is of type D. 
Assume the type of $\sigma$ is $(1, p)$. If $p = 5$ or  not of the form $(r^k - 1)/(r -1)$,  $r$  a prime power,  
then $\Oc$ is sober; otherwise $\Oc$ is of type D. 
 
\smallbreak \noindent $\diamond$ \cite[4.1]{AFGV-ampa} For all other types, $\Oc$ is of type D, hence it collapses.

\smallbreak \noindent $\circ$ \cite{ACG-I} Let $n \ge 2$ and $q$ be a prime power. Let $x\in \PSL_n(q)$ not semisimple and $\Oc = \Oc^{\PSL_n(q)}_{x}$.
The type of a unipotent element are the sizes of its Jordan blocks.

\smallbreak \noindent $\diamond$ Assume $x$ is unipotent. If $x$ is either of type $(2)$ and $q$ is even or not a square, or 
of type $(3)$ and $q =2$, then $\Oc$ is sober. If $x$ is either of type $(2,1)$ and $q$ is even, or  of type $(2,1,1)$ and $q=2$
then $\Oc$ is cthulhu. If $x$ is of type $(2,1, 1)$ and $q > 2$ is even, then $\Oc$ is not of type D, but it is open if it is of type F.
 
\smallbreak \noindent $\diamond$ Otherwise, $\Oc$ is either of type D or of type F, hence it collapses.

\smallbreak \noindent $\circ$ \cite{AFGV-espo, FV} Let $\Oc$ be a conjugacy class in a sporadic simple group $G$.
If $\Oc$ appears in Table \ref{tab:notD}, then $\Oc$ is not of type D. 
If $G = M$ is the Monster and $\Oc$ is one of \textup{32A, 32B, 41A, 46A, 46B, 47A, 47B, 59A, 59B, 69A, 69B, 71A, 71B,
	87A, 87B, 92A, 92B, 94A, 94B}, then it is open  whether $\Oc$ is of type D.
	Otherwise, $\Oc$ is of type D.

\begin{table}[t]
\caption{Classes in sporadic simple groups not of type D}\label{tab:notD}
\begin{center}
\begin{tabular}{|c|c|c|c|}
\hline Group & Classes & Group & Classes\\
\hline $T$ &  \textup{2A} & $Co_{3}$ &  \textup{23A, 23B}\\
\hline $M_{11}$ & \textup{8A, 8B, 11A, 11B} &  $J_{1}$ &  \textup{15A, 15B, 19A, 19B, 19C}\\
\hline $M_{12}$ &  \textup{11A, 11B} & $J_{2}$ &  \textup{2A, 3A}\\
\hline $M_{22}$ &  \textup{11A, 11B} & $J_{3}$ &  \textup{5A, 5B, 19A, 19B}\\
\hline $M_{23}$ &  \textup{23A, 23B} & $J_4$ &  \textup{29A, 43A, 43B, 43C}\\
\hline $M_{24}$ &  \textup{23A, 23B}& $Ly$ &  \textup{37A, 37B, 67A, 67B, 67C}\\
\hline $Ru$ &  \textup{29A, 29B}& $O'N$ &  \textup{31A, 31B}\\
\hline $Suz$ &  \textup{3A} & $Fi_{23}$ & \textup{2A}\\
\hline $HS$ &  \textup{11A, 11B} & $Fi_{22}$ & \textup{2A, 22A, 22B}\\
\hline $McL$ &  \textup{11A, 11B} & $Fi'_{24}$ & \textup{29A, 29B}\\
\hline $Co_{1}$ &  \textup{3A} & $B$ & \textup{2A, 46A, 46B, 47A, 47B}\\
\hline $Co_{2}$ &  \textup{2A, 23A, 23B}  \\
\cline{1-2}
\end{tabular}
\end{center}
\end{table}

\subsection{Generation in degree one}\label{subsec:gen-deg-1}

Here is the scheme of proof proposed in \cite{AS-adv} to attack Conjecture \ref{conj:gendegone}: 
Let $T$ be a finite-dimensional graded  Hopf algebra in $\ydk$  with $T^0 = \kk$ and generated as algebra by $T^1$. We have a commutative diagram of Hopf algebra maps
$\xymatrix@R=2pt@C=10pt{T \ar@{->>}^{\pi}[rr]  & & \toba(V) \\ & T(V)\ar@{->}[ul]^p \ar@{->}[ur]  &}
$. To show that $\pi$ is injective, take a generator $r$ (or a family of generators) of $\J(V)$ such that $r\in \Pc(T(V))$
and consider the Yetter-Drinfeld submodule $U = \kk r\oplus V$ of $T(V)$; if $\dim \toba(U) = \infty$, then $p(r) = 0$. 
Then $p$ factorizes through $T(V)/\J_1(V)$, where $\J_1(V)$ is the ideal generated by primitive generators of $\J(V)$, and so on.

The Conjecture has been verified in all known examples in characteristic 0 (it is false in positive characteristic or for infinite-dimensional Hopf algebras).  
\begin{theorem}        
A finite-dimensional pointed Hopf algebra $H$ is generated by group-like and skew-primitive elements if either of the following holds:

\smallbreak \noindent $\diamond$ \emph{\cite{A-crelle}}  The infinitesimal braiding is of diagonal type, e.~g. $G(H)$ is abelian.
        
\smallbreak \noindent $\diamond$ \emph{\cite{AG-ama, GG}.} The infinitesimal braiding of $H$ is any of $(\Oc^m_2, -1)$, $(\Oc^m_2, \chi)$ ($m= 3,4, 5$), $(X_{4, \omega}, -1)$,
        $(X_{5, 2}, -1)$, $(X_{5, 3}, -1)$, $(X_{7, 3}, -1)$, $(X_{7, 5}, -1)$.
\end{theorem}

\smallbreak\subsection{Liftings}\label{subsec:deformations}
We address here Question \ref{item:method-lifting} in \S \ref{subsec:chevalley}. Let $X$ be a finite rack and $q: X\times X \to \G_{\infty}$
a 2-cocycle. A Hopf algebra $H$ is a \emph{lifting} of $(X, q)$ if $H_0$ is a Hopf subalgebra, $H$ is generated by $H_1$ and 
its infinitesimal braiding is a realization of $(\kk X,c^{q})$. See \cite{GIV} for liftings in the setting of copointed Hopf algebras.

We start  discussing realizations of braided vector spaces as Yetter-Drinfeld modules. Let $\theta\in \N$ and $\I = \I_\theta$.
First, a \emph{YD-datum of diagonal type} is a collection
\begin{align}\label{eq:YD-datum-diagonal}
\D = ((q_{ij})_{i, j\in \I}, G, (g_i)_{i\in \I}, (\chi_i)_{i\in \I}),
\end{align}
where $q_{ij} \in \G_{\infty}$, $q_{ii}\neq 1$, $i, j\in \I$;  $G$ is a finite group; ${g_i} \in Z(G)$; ${\chi_i} \in \widehat{G}$, $i\in \I$; 
such that $q_{ij} := \chi_j(g_i)$, $i, j\in \I$.
Let $(V,c)$ be the braided vector space of diagonal type with matrix  $(q_{ij})$
in the basis $(x_i)_{i\in \I_\theta}$. Then $V\in \ydgg$ by declaring $x_i \in V_{g_i}^{\chi_i}$, $i\in \I$. More generally,
a \emph{YD-datum of rack type} \cite{AG-ama, MS} is a collection
\begin{align}\label{eq:YD-datum-rack}
\D = (X, q, G, \cdot, g, \chi),
\end{align}
where $X$ is a finite rack; $q: X\times X \to \G_{\infty}$ is a $2$-cocycle; $G$ is a finite group;
$\cdot$ is an action of $G$ on $X$; $g: X \to G$ is  equivariant 
with respect to the conjugation in $G$; and $\chi = (\chi_i)_{i\in X}$ is a family  of 1-cocycles
$\chi_i: G \to \kk^{\times}$ (that is, $\chi_i(ht) =\chi_i(t) \chi_{t \cdot i}(h)$, for
all $i\in X$, $h, t\in G$) such that  $g_i \cdot j = i\trid j$ and
$\chi_i(g_j) = q_{ij}$ for all $i, j\in X$. 
Let $(V, c) = (\kk X, c^q)$ be the associated braided vector space. 
Then $V$ becomes an object in $\ydgg$ by $\delta(x_i) = g_i \otimes x_i$ and 
$t\cdot x_{i} = \chi_i(t) x_{t\cdot i}$, $t \in G$, $i\in X$.

Second, let  $\D$ be a YD-datum of either diagonal or rack type and $V \in \ydgg$ as above; let $\T(V) := T(V)\# \kk G$.
The desired liftings are quotients of $\T(V)$; write $a_i$ in these quotients
instead of $x_i$ to distinguish them from the elements in $\toba(V)\# \kk G$. 
Let $\cG$ be a minimal set  of  generators of $\J(V)$, assumed homogeneous both for the $\N$- and the $G$-grading.
Roughly speaking, the deformations will be defined by replacing the relations $r = 0$ by $r = \phi_r$, $r\in \cG$,
where $\phi_r\in \T(V)$ belongs to a lower term of the coradical filtration, and the ideal $\J_{\phi}(V)$ generated by $\phi_r$, $r\in \cG$,
is a Hopf ideal. The problem is  to describe the $\phi_r$'s and to check that $\T(V)/\J_{\phi}(V)$ has the right dimension. 
If $r\in \Pc(T(V))$  has $G$-degree $g$, then $\phi_r = \lambda(1 - g)$ for some $\lambda\in \kk$; depending on the action of $G$ on $r$,
it may happen that $\lambda$ should be 0. In some cases, all $r\in \cG$ are primitive, so all deformations can be described; see \cite{AS-jalg} 
for quantum linear spaces (their liftings can also be presented as Ore extensions \cite{BDG}) 
and the Examples \ref{exa:lifting-transpositions-s3} and \ref{exa:lifting-transpositions-s4}. 
But in most cases, not all $r\in \cG$ are primitive and some recursive construction of the deformations is needed. This was achieved in
\cite{AS-cambr} for diagonal braidings of Cartan type $A_n$, with explicit formulae, and in \cite{AS-ann} for diagonal braidings of finite Cartan type,
with recursive formulae. Later it was observed that the so obtained liftings are cocycle deformations of $\toba(V)\# \kk G$, see e.g. \cite{M4}.
This led to the strategy in \cite{AAGMV}:   pick an adapted stratification $\cG =  \cG_0 \cup \cG_1 \cup\dots \cup \cG_{N}$ \cite[5.1]{AAGMV};
then construct recursively the deformations of $T(V) / \langle \cG_0 \cup \cG_1 \cup\dots \cup \cG_{k-1}\rangle$ 
by determining the cleft extensions of the deformations in the previous step and  applying the theory of Hopf bi-Galois extensions \cite{S-bigal}. 
In the Examples below, $\chi_i = \chi \in \widehat{G}$  for all $i\in X$ by \cite[3.3 (d)]{AG-ama}.  

\begin{example}\label{exa:lifting-transpositions-s3} \cite{AG-ama, GIV}
Let $\D = (\Oc^3_2, -1, G, \cdot, g, \chi)$ be a YD-datum. Let $\lambda \in \kk^2$ be such that
\begin{align}
\label{norm3} \lambda_1 &= \lambda_2 = 0,&  &\text{ if }\chi^2\neq \varepsilon;
\\
\label{norm1}
\lambda_1 &= 0, &  &\text{ if } g_{12}^2 = 1;  &
\lambda_2 &= 0, & &\text{ if } g_{12}g_{13} = 1.
\end{align}
Let $\A = \A(\D, \lambda)$ be the quotient of $\T(V)$ by the relations
\begin{align}
\label{relsym31}
a_{12}^2 &= \lambda_1(1-g_{12}^2), 
\\\label{relsym33}
 a_{12}a_{13} + a_{23}a_{12}+ a_{13}a_{23} &= \lambda_2 (1-g_{12} g_{13}).
\end{align}
Then $\A$ is a pointed Hopf algebra, a cocycle deformation of $\gr \A\simeq \toba(V)\#\kk G$ and $\dim \A = 12 |G|$;
$\A(\D, \lambda) \simeq \A(\D, \lambda')$ iff $\lambda = c\lambda'$ for some $c\in \kk^{\times}$.
Conversely, any lifting of $(\Oc^3_2,-1)$ is isomorphic to $\A(\D, \lambda)$ for some YD-datum
  $\D = (\Oc^3_2, -1, G, \cdot, g, \chi)$ and $\lambda \in \kk^2$  satisfying \eqref{norm3}, \eqref{norm1}. 
\end{example}

\begin{example}\label{exa:lifting-transpositions-s4} \cite{AG-ama}
Let $\D = (\Oc^4_2, -1, G, \cdot, g, \chi)$ be a YD-datum.
Let $\lambda \in \kk^3$ be such that
\begin{align}\label{norm3bis}
\lambda_i &= 0,& &i\in \I_3,&  \text{ if }\chi^2&\neq \varepsilon;
\\\label{norm14}
\lambda_1 &= 0, &  &\text{ if } g_{12}^2 = 1;  &
\lambda_2 &= 0, & \text{ if } g_{12}g_{34} &= 1; &
\lambda_3 &= 0, & \text{ if } g_{12}g_{13} &= 1.
\end{align}
Let $\A = \A(\D, \lambda)$ be the quotient of $\T(V)$ by the relations
\begin{align}
\label{relsym41}
a_{12}^2 &= \lambda_1(1-g_{12}^2),
\\\label{relsym42}
a_{12} a_{34} +a_{34} a_{12}  &= \lambda_2 (1-g_{12} g_{34}), 
\\\label{relsym43}
a_{12}a_{13} + a_{23}a_{12}+ a_{13}a_{23} &= \lambda_3 (1- g_{12} g_{13}).
\end{align}
Then $\A$ is a pointed Hopf algebra, a cocycle deformation of  $\toba(V)\#\kk G$ and $\dim \A = 576 |G|$;
$\A(\D, \lambda) \simeq \A(\D, \lambda')$ iff $\lambda = c\lambda'$ for some $c\in \kk^{\times}$.
Conversely, any lifting of $(\Oc^4_2, -1)$ is isomorphic to $\A(\D, \lambda)$ for some YD-datum 
$\D = (\Oc^4_2, -1, G, \cdot, g, \chi)$ and $\lambda \in \kk^3$  satisfying \eqref{norm3bis}, \eqref{norm14}. 
\end{example}

\begin{example}\label{exa:lifting-affine-f4} \cite{GIV}
Let $\D = (X_{4,\omega}, -1, G, \cdot, g, \chi)$ be a YD-datum. Let $\lambda \in \kk^3$ be such that
\begin{align}
\label{eqn:cond1-A} \lambda_1 &= \lambda_2 = 0,&   \text{ if }\chi^2&\neq \varepsilon;&  \lambda_3 &=  0,  \text{ if }\chi^6\neq \varepsilon;\\
\label{eqn:cond1-L} \lambda_1 &=0, \text{ if }\,  g_0^2=1, & \lambda_2 &=0, \text{ if } \, g_0g_1 =1, & \lambda_3&=0,  \text{ if }\, g_0^3g_1^3=1.
\end{align}
Let $\A = \A(\D, \lambda)$ be the quotient of $\T(V)$ by the relations

\begin{align}
x_0^2 &= \lambda_1(1- g_0^2), \\ 
x_0x_1 + x_{\omega} x_0 + x_1x_{\omega}  &= \lambda_2(1 -g_0g_1) \\
(x_{\omega}x_1x_0)^2 + (x_1x_0x_{\omega})^2 + (x_0x_{\omega}x_1)^2  &= \zeta_6 - \lambda_3(1-g_0^3g_1^3),& &\text{where}
\end{align} 
\begin{align*}&\zeta_6=\lambda_2(x_{\omega}x_1x_0x_{\omega}+ x_1x_0x_{\omega}x_1
+x_0x_{\omega}x_1x_0)-\lambda_2^3(g_0g_1-g_0^3g_1^3)
\\\notag
+ & \lambda_1^2g_0^2\big(g_{1 + \omega}^2(x_{\omega}x_3 +x_0x_{\omega}) +g_1g_{1 + \omega}(x_{\omega}x_1+
x_1x_3) +g_1^2(x_1x_0+x_0x_3)\big)
\\\notag
- & 2\lambda_1^2g_0^2(x_0x_3- x_{\omega}x_3   -x_1x_{\omega} + x_1x_0) - 2\lambda_1^2g_{\omega}^2(x_{\omega}x_3
-x_1x_3+x_0x_{\omega} - x_0x_1 )
\\\notag
- & 2\lambda_1^2g_1^2(x_{\omega}x_1 +  x_1x_3 +  x_1x_{\omega} - x_0x_3 + x_0x_1) 
\\\notag
+ & \lambda_2\lambda_1(g_{\omega}^2x_0x_3+ g_1^2x_{\omega}x_3+ g_0^2x_1x_3)
+\lambda_2^2g_0g_1(x_{\omega}x_1 + x_1x_0 + x_0x_{\omega}-\lambda_1)
\\ \notag
- & \lambda_2\lambda_1^2(3g_0^3g_{1 + \omega}-2 g_0g_1^3- g_0^2g_{\omega}^ 2
-2 g_0^3g_1+g_{\omega}^2 - g_1^2+g_0^2)
\\\notag
 - &\lambda_2(\lambda_1-\lambda_2)\big(\lambda_1\,g_0^2(g_{1 + \omega}^2+g_1g_{1 + \omega} +g_1^2+2g_0g_1^3)
+x_{\omega}x_1 + x_1x_0 + x_0x_{\omega}\big).
\end{align*}
Then $\A$ is a pointed Hopf algebra, a cocycle deformation of $\gr \A\simeq \toba(V)\#\kk G$ and $\dim \A = 72 |G|$;
$\A(\D, \lambda) \simeq \A(\D, \lambda')$ iff $\lambda = c\lambda'$ for some $c\in \kk^{\times}$.
Conversely, any lifting of $(X_{4,\omega}, -1)$ is isomorphic to $\A(\D, \lambda)$ for some  YD-datum 
$\D = (X_{4,\omega}, -1, G, \cdot, g, \chi)$ and $\lambda \in \kk^3$  satisfying \eqref{eqn:cond1-A}, \eqref{eqn:cond1-L}. 
\end{example}

\begin{example}\label{exa:lifting-affine-f5} \cite{GIV}
Let $\D = (X_{5,2}, -1, G, \cdot, g, \chi)$ be a YD-datum. Let $\lambda \in \kk^3$ be such that
\begin{align}
\label{eqn:cond1-Aa} \lambda_1 &= \lambda_2 = 0,&   \text{ if }\chi^2&\neq \varepsilon;&  \lambda_3 &=  0,  \text{ if }\chi^4\neq \varepsilon;\\
\label{eqn:cond1-La} \lambda_1 &=0, \text{ if }\,  g_0^2=1, & \lambda_2 &=0, \text{ if } \, g_0g_1 =1, & \lambda_3&=0,  \text{ if }\,g_0^2g_1g_2 = 1.
\end{align}
Let $\A = \A(\D, \lambda)$ be the quotient of $\T(V)$ by the relations
\begin{align}
&x_0^2 = \lambda_1(1- g_0^2), \\ 
&x_0x_1 + x_2x_0 + x_3x_2 + x_1x_3  = \lambda_2(1 -g_0g_1), \\
&(x_1x_0)^2 +  (x_0x_1)^2  = \zeta_4 - \lambda_3(1-g_0^2g_1g_2),& &\text{where}
\end{align}
$\zeta_4= \lambda_2\, (x_1x_0 +
x_0x_1)+\lambda_1\, g_1^2(x_3x_0+ x_2x_3) - \lambda_1\, g_0^2(x_2x_4+ x_1x_2) +
\lambda_2\lambda_1\,g_0^2(1- g_1g_2)$.
Then $\A$ is a pointed Hopf algebra, a cocycle deformation of $\gr \A\simeq \toba(V)\#\kk G$ and $\dim \A = 1280 |G|$;
$\A(\D, \lambda) \simeq \A(\D, \lambda')$ iff $\lambda = c\lambda'$ for some $c\in \kk^{\times}$.
Conversely, any lifting of $(X_{5,2}, -1)$ is isomorphic to $\A(\D, \lambda)$ for some  YD-datum 
$\D = (X_{5,2}, -1, G, \cdot, g, \chi)$ and $\lambda \in \kk^3$  satisfying \eqref{eqn:cond1-Aa}, \eqref{eqn:cond1-La}.  
\end{example}

\begin{example}\label{exa:lifting-affine-f5-3} \cite{GIV}
Let $\D = (X_{5,3}, -1, G, \cdot, g, \chi)$ be a YD-datum. Let $\lambda \in \kk^3$ be such that
\begin{align}
\label{eqn:cond1-Ab} \lambda_1 &= \lambda_2 = 0,&   \text{ if }\chi^2&\neq \varepsilon;&  \lambda_3 &=  0,  \text{ if }\chi^4\neq \varepsilon;\\
\label{eqn:cond1-Lb} \lambda_1 &=0, \text{ if }\,  g_0^2=1, & \lambda_2 &=0, \text{ if } \, g_1g_0 =1, & \lambda_3&=0,  \text{ if }\,g_0^2g_1g_3 = 1.
\end{align}
Let $\A = \A(\D, \lambda)$ be the quotient of $\T(V)$ by the relations
\begin{align}
&x_0^2 = \lambda_1(1- g_0^2), \\ 
&x_1x_0 + x_0x_2 + x_2x_3 + x_3x_1  = \lambda_2(1 -g_1g_0) \\
&x_0x_2x_3x_1 + x_1x_4x_3x_0 = \zeta'_4 - \lambda_3(1-g_0^2g_1g_3),
\end{align}
where  
$\zeta'_4= \lambda_2\, (x_0x_1 + x_1x_0) -\lambda_1\, g_1^2(x_3x_2 +x_0x_3) - \lambda_1\, 
g_0^2(x_3x_4 +x_1x_3)  + \lambda_1\lambda_2( g_1^2 + 
g_0^2- 2 g_0^2g_1g_3)$.
Then $\A$ is a pointed Hopf algebra, a cocycle deformation of $\gr \A\simeq \toba(V)\#\kk G$ and $\dim \A = 1280 |G|$;
$\A(\D, \lambda) \simeq \A(\D, \lambda')$ iff $\lambda = c\lambda'$ for some $c\in \kk^{\times}$.
Conversely, any lifting of $(X_{5,3}, -1)$ is isomorphic to $\A(\D, \lambda)$ for some  YD-datum 
$\D = (X_{5,3}, -1, G, \cdot, g, \chi)$ and $\lambda \in \kk^3$  satisfying \eqref{eqn:cond1-Ab}, \eqref{eqn:cond1-Lb}.   
\end{example}

\section{Pointed Hopf algebras}\label{sec:pointed}

\smallbreak\subsection{Pointed Hopf algebras with abelian group}\label{subsec:pointed-ab}
Here is a classification from \cite{AS-ann}.
Let $\D = ((q_{ij})_{i, j\in \I_{\theta}}, \Gamma, (g_i)_{i\in \I_\theta}, (\chi_i)_{i\in \I_{\theta}})$ be a YD-datum of diagonal type
as in  \eqref{eq:YD-datum-diagonal} with $\Gamma$ a finite abelian group and let $V\in \ydg$ be the corresponding realization. 
We say that $\D$ is a \emph{Cartan datum} if 
there is a Cartan matrix (of finite type) $\ba = (a_{ij})_{i, j\in \I_{\theta}}$  
such that $q_{ij} q_{ji} = q_{ii}^{a_{ij}}$, $i\neq j \in \I_{\theta}$.

\smallbreak
Let $\Phi$ be the root system associated to  $\ba$, $\alpha_1,\dots,\alpha_{\theta}$  a choice of simple roots, 
$\mathcal{X}$ the set of connected components of the Dynkin diagram of  $\Phi$ and set 
$i \sim j$ whenever $\alpha_i, \alpha_j$ belong to the same  $J\in \mathcal X$. We consider two classes of parameters:
\begin{itemize}\renewcommand{\labelitemi}{$\circ$}\itemsep1pt \parskip0pt \parsep0pt
 \item $ \lambda =(\lambda_{ij})_{\substack{i < j \in \I_\theta,\\ i\not\sim j}}$ is a family in  $\{0,1\}$ with $\lambda_{ij} =0$ when
$g_ig_j =1$ or $\chi_i \chi_j \neq \varepsilon$. 

 \item 
$\mu=(\mu_{\alpha})_{\alpha \in \Phi^+}$ is a family in  $\kk$
with $\mu_{\alpha} =0$ when $\supp \alpha \subset J$, $J \in \cX$,  and  $g_{\alpha}^{N_J} =1$ or $\chi_{\alpha}^{N_J} \neq\varepsilon$. 
Here $N_J =\ord q_{ii}$ for an arbitrary $i\in J$.
\end{itemize}
We attach a family ($u_{\alpha}(\mu))_{\alpha \in \Phi^+}$ in $\kk\Gamma$ to the parameter $\mu$ , defined 
recursively on the length of $\alpha$, starting by $u_{\alpha_i}(\mu) = \mu_{\alpha_i}(1 - g_i^{N_i})$.
From all these data we define a Hopf algebra   $\ug(\D,\lambda,\mu)$ as the quotient of 
$\T(V) = T(V) \#\kk \Gamma$ by  the relations
\begin{align}
ga_i g^{-1} &= \chi_i(g)a_i, & & \\
\ad_c(a_i)^{1 - a_{ij}}(a_j) &= 0,&  &i \neq j, i\sim j,\\
\ad_c(a_i)(a_j) &= \lambda_{ij}(1 - g_ig_j),& & i<j ,i \nsim j,\\
 a_{\alpha}^{N_J} &= u_{\alpha}(\mu).
\end{align}

\begin{theorem} The Hopf algebra ${\ug(\D,\lambda,\mu)}$ is pointed, $G({\ug(\D,\lambda,\mu)}) \simeq \Gamma$ and 
$\dim {\ug(\D,\lambda,\mu)} = \prod_{J \in \cX} N_J^{|\Phi_J^+|} |\Gamma|$.
Let $H$ be a pointed finite-dimensional Hopf algebra and set $\Gamma = G(H)$. 
Assume that the prime divisors of $\vert\Gamma\vert$ are $>7$. Then there exists a Cartan datum
$\D$ and parameters $\lambda$ and $\mu$ such that $H \simeq \ug(\D,\lambda,\mu)$.
It is known when two Hopf algebras ${\ug(\D,\lambda,\mu)}$ and ${\ug(\D',\lambda',\mu')}$ are isomorphic.
 \end{theorem}

The proof offered in \cite{AS-ann} relies on \cite{AS-adv, H-inv, L1, L2}. 
Some comments: the hypothesis on $\vert \Gamma\vert$ forces the infinitesimal braiding $V$ of $H$ to be of Cartan type, 
and the relations of $\toba(V)$ to be just quantum Serre and powers of root vectors. The quantum Serre relations
are not deformed in the liftings, except those linking different components of the Dynkin diagram; the powers of the root vectors
are deformed to the $u_{\alpha}(\mu)$ that belong to the coradical. All this can fail without the hypothesis, see \cite{helbig} for examples
in rank 2.

\smallbreak\subsection{Pointed Hopf algebras with non-abelian group}\label{subsec:pointed-non-ab}
We present some classification results of pointed Hopf algebras with non-abelian group.
We say that a finite group $G$ \emph{collapses} whenever any finite-dimensional pointed Hopf algebra $H$ with $G(H)\simeq G$ is isomorphic to $\kk G$.

\smallbreak \noindent $\bullet$  \cite{AFGV-ampa, AFGV-espo}  Let $G$ be either  $\am$, $m\geq 5$, or a
sporadic simple group, different from  $Fi_{22}$, the Baby Monster $B$ or the Monster $M$.
Then $G$ collapses.
 
The proof uses  \S \ref{subsubsec:nichols-rack-collapsing};  the remaining Yetter-Drinfeld modules are discarded 
considering abelian subracks of the supporting conjugacy class and  the list in \cite{H-classif}.

\smallbreak \noindent $\bullet$ \cite{AHS} Let $V = M(\Oc^3_2, \sgn)$ and let  $\D$ be the corresponding YD-datum.
Let $H$ be a finite-dimensional pointed Hopf algebra  with $G(H)\simeq \st$.
Then $H$ is isomorphic either to $\kk\st$, or to $\A(\D, 0) = \toba(V)\# \kk\st$, or to $\A(\D, (0,1))$, cf. Example \ref{exa:lifting-transpositions-s3}.

\smallbreak \noindent $\bullet$ \cite{GG} Let $H\not\simeq \kk\sk$ be a finite-dimensional pointed Hopf algebra  with $G(H)\simeq \sk$. 
Let $V_1 = M(\Oc^4_2, \sgn \otimes \id)$, $V_2 = M(\Oc^4_2, \sgn \otimes \sgn)$, $W = M(\Oc^4_4, \sgn \otimes \id)$, with corresponding data $\D_1$, 
$\D_2$ and $\D_3$.
Then $H$ is isomorphic to one of
\begin{align*}
&\A(\D_1, (0,\mu)),\quad \mu\in \kk^2;&
&\A(\D_2, t), \quad t\in \{0,1\};& &\A(\D_3, \lambda),\quad  \lambda\in\kk^2.
\end{align*}
Here $\A(\D_1, (0,\mu))$ is as in Example \ref{exa:lifting-transpositions-s4}; 
$\A(\D_2, t)$ is the quotient of $\T(V_2)$ by the relations
$ a_{12}^2 = 0$, $a_{12}a_{34} - a_{34}a_{12} = 0$, 
$a_{12}a_{23} - a_{13}a_{12} - a_{23}a_{13} = t(1- g_{(12)}g_{(23)})$
and
$\A(\D_3, \lambda)$ is the quotient of $\T(W)$ by the relations
\begin{align*}
& a_{(1234)}^2= \lambda_1(1-g_{(13)}g_{(24)}); \qquad a_{(1234)}a_{(1432)}+a_{(1432)}a_{(1234)}=0;\\
&a_{(1234)}a_{(1243)}+a_{(1243)}a_{(1423)}+a_{(1423)}a_{(1234)}=\lambda_2(1-g_{(12)}g_{(13)}).
\end{align*}
Clearly $\A(\D_1, 0) = \toba(V_1)\# \kk\sk$, $\A(\D_2, 0) = \toba(V_2)\# \kk\sk$, $\A(\D_3, 0) = \toba(W)\# \kk\sk$. 
Also $\A(\D_1, (0,\mu)) \simeq \A(\D_1, (0,\nu))$ iff $\mu = c\nu$ for some $c \in \kk^{\times}$, and $\A(\D_3, \lambda) 
\simeq \A(\D_3, \kappa)$ iff $\lambda = c\kappa$ for some $c \in \kk^{\times}$.

\smallbreak \noindent $\bullet$ \cite{AFGV-ampa, GG} Let $H$ be a finite-dimensional pointed Hopf algebra  with $G(H)\simeq \sco$, but $H\not\simeq \kk\sco$.
It is not known whether $\dim \toba(\Oc_{2,3}^5,\sgn\otimes\varepsilon) <\infty$. 
Let $\D_1$, $\D_2$ be the data corresponding to $V_1 = M(\Oc^5_2, \sgn \otimes \id)$, $V_2 = M(\Oc^5_2, \sgn \otimes \sgn)$. 
If the infinitesimal braiding of 
$H$ is not $M(\Oc_{2,3}^5,\sgn\otimes\varepsilon)$, then $H$ is isomorphic to  one of
$\A(\D_1, (0,\mu))$, $\mu\in \kk^2$ (defined as in  Example \ref{exa:lifting-transpositions-s4}),
or $\toba(V_2)\# \kk\sco$, or $\A(\D_2, 1)$ (defined as above).

\smallbreak \noindent $\bullet$ \cite{AFGV-ampa} Let $m > 6$.
Let $H\not\simeq \kk\sm$ be a finite-dimensional pointed Hopf algebra  with $G(H)\simeq \sm$.
Then the infinitesimal braiding of $H$ is  $V = M(\Oc, \rho)$, where the type of $\sigma$ is $(1^{m-2}, 2)$ and $\rho = \rho_1 \otimes \sgn$, 
$\rho_1 =\sgn$ or $\varepsilon$; it is an open question whether $\dim \toba(V) < \infty$, see Example \ref{exa:nichols-sim}. If $m = 6$, there 
are two more Nichols algebras with unknown dimension corresponding to the class of type $(2^3)$,  but they are conjugated to those
of type $(1^4, 2)$ by the outer automorphism of $\sei$.
    
\smallbreak \noindent $\bullet$ \cite{FG}  Let $m \ge 12$, $m = 4h$ with $h\in \N$.
Let $G = \dm$ be the dihedral group of order $2m$. Then there are infinitely many finite-dimensional Nichols algebras
in $\ydgg$; all of them are exterior algebras as braided Hopf algebras.
Let $H$ be a finite-dimensional pointed Hopf algebra  with $G(H)\simeq \dm$, but $H\not\simeq \kk\dm$.
Then $H$ is a lifting of an exterior algebra, and there infinitely many such liftings.

\smallbreak\subsubsection{Copointed Hopf algebras}\label{subsubsec:co-pointed}
We say that a semisimple Hopf algebra $K$ \emph{collapses} if any finite-dimensional  Hopf algebra $H$ with $H_0\simeq K$ is isomorphic to $K$.
Thus, if $G$ collapses, then $\kk^G$ and $(\kk G)^F$ collapse, for any twist $F$.
Next we state the classification of the finite-dimensional copointed Hopf algebras over $\st$  \cite{AV}.
Let $V = M(\Oc^3_2, \sgn)$ as a Yetter-Drinfeld module over $\kk^{\st}$.
Let $\lambda \in\kk^{\Oc_2^3}$ be such that $\displaystyle\sum_{(ij)\in\Oc_2^3}\lambda_{ij}=0$. 
Let $\vg = \vg(V, \lambda)$ be the quotient of $T(V)\#\kk^{\st}$ by the relations
$\xij{13}\xij{23}+\xij{12}\xij{13} + \xij{23}\xij{12} = 0$, $\xij{23}\xij{13}+\xij{13}\xij{12}+\xij{12}\xij{23} =0$, 
\newline $\xij{ij}^2 = \sum_{g\in\st}(\lambda_{ij} - \lambda_{g^{-1}(ij)g})\delta_g$,  for $(ij)\in \Oc_2^3$.
Then $\vg$ is a Hopf algebra of dimension 72 and $\gr \vg \simeq \toba(V)\# \kk^{\st}$. 
Any finite-dimensional copointed Hopf algebra $H$ with $H_0 \simeq \kk^{\st}$ is isomorphic to
$\vg(V, \lambda)$ for some $\lambda$ as above; $\vg(V, \lambda)\simeq \vg(V, \lambda')$ iff $\lambda$
and $\lambda'$ are conjugated under $\kk^{\times}\times\Aut \st$.


\begin{thebibliography}{80}


\bibitem{bariloche} Andruskiewitsch, N., About finite dimensional Hopf algebras, \emph{Contemp. Math.}
\textbf{294} (2002), 1--54.

\bibitem{AA} Andruskiewitsch, N.,   Angiono, I.,  Weyl groupoids, contragredient Lie superalgebras and Nichols algebras,
in preparation.

\bibitem{AAGMV} Andruskiewitsch, N.,   Angiono, I., Garc\'{\i}a Iglesias, A.,  Masuoka, A., Vay, C., Lifting via cocycle deformation,
\emph{J. Pure Appl. Alg.}  \textbf{218} (2014), 684--703.


\bibitem{ACG-I}  Andruskiewitsch, N.,  Carnovale, G., Garc\'ia, G. A.,
Finite-dimensional pointed Hopf algebras over finite simple groups of Lie type  I, \texttt{arXiv:1312.6238}.

\bibitem{AC} Andruskiewitsch, N., Cuadra, J., On the structure of (co-Frobenius) Hopf algebras,
\emph{J. Noncommut. Geom.}  \textbf{7} (2013),   83--104.

\bibitem{AEG} Andruskiewitsch, N., Etingof, P., Gelaki, S.,
Triangular Hopf algebras with the Chevalley property, \emph{Michigan Math. J.} \textbf{49} (2001), 277--298.

\bibitem{AFGV-ampa} Andruskiewitsch, N., Fantino, F.,  Gra\~na, M.,   Vendramin, L., Finite-di\-mensional pointed
Hopf algebras with alternating groups are trivial, \emph{Ann. Mat. Pura Appl. (4)},  \textbf{190}  (2011), 225--245.

\bibitem{AFGV-espo} Andruskiewitsch, N., Fantino, F.,  Gra\~na, M.,   Vendramin, L., Pointed Hopf algebras over the sporadic simple groups,
\emph{J. Algebra} \textbf{325} (2011),  305--320.

\bibitem{AG-compo}  Andruskiewitsch, N., Garc\'\i a, G. A., Finite subgroups of a simple quantum group,
\emph{Compositio Math.} \textbf{145} (2009), 476--500.

\bibitem{AG-adv} Andruskiewitsch, N.,   Gra\~na, M.,
From racks to pointed Hopf algebras, \emph{Adv. Math.}  \textbf{178}  (2003), 177--243.

\bibitem{AG-ama} Andruskiewitsch, N.,   Gra\~na, M., Examples of
liftings of Nichols algebras over racks, \emph{AMA Algebra Montp. Announc. (electronic)}, Paper {\bf 1},  (2003).

\bibitem{AHS} Andruskiewitsch, N., Heckenberger, I.,   Schneider, H.-J., 
The Nichols algebra of a semisimple Yetter-Drinfeld module, \emph{Amer. J. Math.} \textbf{132} (2010),  1493--1547.

\bibitem{AS-jalg} Andruskiewitsch, N.,   Schneider, H.-J.,  Lifting of quantum linear spaces and pointed Hopf algebras of order $ p^3$, 
\emph{J. Algebra} \textbf{209} (1998), 658--691.

\bibitem{AS-adv} Andruskiewitsch, N.,   Schneider, H.-J.,  Finite quantum groups
and Cartan matrices, \emph{Adv. Math.} \textbf{154} (2000), 1--45.

\bibitem{AS-cambr} Andruskiewitsch, N.,   Schneider, H.-J.,  Pointed Hopf algebras, in  \emph{New directions in Hopf algebras},
MSRI series, Cambridge Univ. Press (2002), 1--68.


\bibitem{AS-ann} Andruskiewitsch, N.,   Schneider, H.-J., On the classification of finite-di\-mensional pointed Hopf
algebras, \emph{Ann. of Math.} \textbf{171} (2010),  375--417.


\bibitem{AV} Andruskiewitsch, N.,   Vay, C., Finite dimensional Hopf algebras over the dual group algebra of the symmetric group in 3 letters,
\emph{Comm. Algebra} \textbf{39} (2011), 4507--4517.

\bibitem{A-jems} Angiono, I.,  A presentation by generators and relations of Nichols algebras 
of diagonal type and convex orders on root systems,  \emph{J. Europ. Math. Soc.}, to appear.

\bibitem{A-crelle}  Angiono, I., On Nichols algebras of diagonal type, \emph{J. Reine Angew. Math.}  \textbf{683} (2013), 189--251.

\bibitem{BDG}
Beattie, M., D\u{a}sc\u{a}lescu, S., Gr\"unenfelder, L.,  On the number of types of finite-dimensional Hopf algebras,
\emph{Invent. Math.} {\bf 136} (1999), 1--7.

\bibitem{BG}
Beattie, M.,  Garc\'\i a, G. A., Classifying Hopf algebras of a given dimension,
\emph{Contemp. Math.} {\bf 585} (2013), 125--152. 

\bibitem{CH-at most 3} Cuntz, M.,   Heckenberger, I., Weyl groupoids with at most three objects, \emph{J. Pure Appl. Algebra} \textbf{213} (2009), 1112--1128.

\bibitem{CH} Cuntz, M.,   Heckenberger, I., Finite Weyl groupoids.  \emph{J. Reine Angew. Math.}, to appear.

\bibitem{DT} Doi, Y., Takeuchi, M., Multiplication alteration by two-cocycles. The 
quantum version, \emph{Comm. Algebra} \textbf{22} (1994), 5715--5732.  

\bibitem{Dr1} Drinfeld, V. G., Quantum groups, 
\emph{Proc. Int. Congr. Math., Berkeley/Calif. 1986}, Vol. 1, 798--820 (1987).

\bibitem{Dr-quasi} Drinfeld, V. G., Quasi-Hopf  algebras, 
\emph{Leningrad Math.  J.} \textbf{1} (1990), 1419--1457.


\bibitem{EG-frob-prop}  Etingof, P., Gelaki, S.,  Some properties of finite-dimensional semisimple Hopf
algebras, \emph{Math. Res. Lett.} \textbf{5} (1998), 191--197.

\bibitem{EG-clas-triang-ss}  Etingof, P., Gelaki, S.,
The classification of triangular semisimple and co\-semisimple Hopf algebras over an algebraically closed field,
\emph{Int. Math. Res. Not. } \textbf{2000, No. 5} (2000), 223--234.

\bibitem{EG-clas-triang}  Etingof, P., Gelaki, S.,  
The classification of finite-dimensional triangular Hopf algebras over an algebraically closed field of characteristic 0, 
\emph{Mosc. Math. J.} \textbf{3} (2003), 37--43.

\bibitem{EGNO}  Etingof, P.,  Gelaki, S.,  Nikshych, D., Ostrik, V., \emph{Tensor categories}, book to appear.


\bibitem{EN}  Etingof, P.,  Nikshych, D., Dynamical quantum groups at roots of 1, \emph{Duke Math. J.}
 \textbf{108} (2001), 135--168.

\bibitem{ENO}  Etingof, P.,  Nikshych, D., Ostrik, V., Weakly group-theoretical and solvable fusion categories, 
\emph{Adv. Math.} \textbf{226} (2011), 176--205.



\bibitem{F} Fantino, F., Conjugacy classes of p-cycles of type D in alternating groups, \emph{Comm. Algebra}, to appear.

\bibitem{FG}  Fantino, F.,  Garc\'ia, G. A., On pointed Hopf algebras over dihedral groups, \emph{Pacific J. Math.} {\bf 252} (2011), 69--91.
 
\bibitem{FV} Fantino, F.,  Vendramin, L., On twisted conjugacy classes of type D in sporadic simple groups,
\emph{Contemp. Math.} \textbf{585} (2013), 247--259.

\bibitem{FK} Fomin, S., Kirillov, K. N.,  Quadratic algebras, Dunkl elements, and Schubert calculus, \emph{Progr. Math.}  \textbf{172} (1999),  146--182.


\bibitem{GN} Galindo, C., Natale, S., Simple Hopf algebras and deformations of finite groups, \emph{Math. Res. Lett.} \textbf{14} (2007), 943--954.


\bibitem{GG} Garc\'{\i}a, G. A., Garc\'{\i}a Iglesias, A.,  Finite dimensional pointed Hopf algebras over $\sk$,
\emph{Israel J. Math.} \textbf{183} (2011), 417--444. 

\bibitem{GIV} Garc\'{\i}a Iglesias, A.,  Vay, C., Finite-dimensional pointed or copointed Hopf
algebras over affine racks, \emph{J. Algebra} \textbf{397} (2014), 379--406.

\bibitem{GNN} Gelaki, S., Naidu, D., Nikshych, D., Centers of graded fusion categories,
\emph{Algebra Number Theory} \textbf{3} (2009), 959--990.

\bibitem{grania1}Gra\~na, M., On Nichols algebras of low dimension, \emph{Contemp. Math.} {\bf 267} (2000), 111--134.

\bibitem{HeckGran07} Gra{\~n}a, M., Heckenberger, I., On a factorization of graded {H}opf
  algebras using {L}yndon words, \emph{J. Algebra} \textbf{314} (2007), 324--343.

\bibitem{H-inv}  Heckenberger, I., The Weyl groupoid of a Nichols algebra of diagonal type, \emph{Invent. Math.} \textbf{164} (2006), 175--188.

\bibitem{H-classif}  Heckenberger, I., Classification of arithmetic root systems, \emph{Adv. Math.} \textbf{220} (2009), 59--124.

\bibitem{HLV}  Heckenberger, I., Lochmann, A.,  Vendramin, L.,
Braided racks, Hurwitz actions and Nichols algebras with many cubic relations, \emph{Transform. Groups} \textbf{17} (2012), 157--194.


\bibitem{HeckSch}  Heckenberger, I., Schneider, H.-J.,  Root systems and Weyl groupoids for Nichols algebras, 
\emph{Proc. Lond. Math. Soc.} \textbf{101} (2010), 623--654.

\bibitem{HS2}  Heckenberger, I., Schneider, H.-J.,  Nichols algebras over groups with finite root system of rank two {I}, 
\emph{ J. Algebra} \textbf{324} (2010), 3090--3114.

\bibitem{HS-israel}  Heckenberger, I., Schneider, H.-J.,  
Right coideal subalgebras of Nichols algebras and the Duflo order on the Weyl groupoid, 
\emph{Israel J. Math.} \textbf{197} (2013), 139--187.

\bibitem{HV1}  Heckenberger, I., Vendramin, L.,
New examples of {N}ichols algebras with finite root system of rank two, \texttt{arXiv:1309.4634}.

\bibitem{HV2}  Heckenberger, I., Vendramin, L.,
The classification of Nichols algebras with finite root system of rank two, \texttt{arXiv:1311.2881}.



\bibitem{HY}  Heckenberger, I., Yamane, H., A generalization of Coxeter groups, root systems, and Matsumoto's theorem, \emph{Math. Z.}
\textbf{259} (2008), 255--276.

\bibitem{helbig}  Helbig, M., On the lifting of Nichols algebras, \emph{Comm. Algebra}
\textbf{40} (2012), 3317--3351.


\bibitem{J} Jimbo, M., A $q$-difference analogue of $U({\mathfrak g})$ 
and the Yang-Baxter equation, \emph{Lett. Math. Phys.} \textbf{10} (1985),  63--69.

\bibitem{JS} Joyal, A., Street, R., Braided tensor categories, \emph{Adv. Math.} \textbf{102} (1993), 20--78.

\bibitem{jo} Joyce, D., Simple quandles, \emph{J. Algebra} \textbf{79} (1982), 307--318.

\bibitem{K-ext} Kac, G.,  Extensions of groups to ring groups, \emph{Math. USSR. Sb.} \textbf{5} (1968), 451--474.

\bibitem{K-libro} Kac, V., \emph{Infinite-dimensional Lie algebras}. 3rd edition. Cambridge Univ. Press,  1990. 

\bibitem{Khar99} Kharchenko, V., A quantum analogue of the  {P}oincar{\'e}--Birkhoff--{W}itt theorem, \emph{Algebra and Logic} \textbf{38} (1999),  259--276.

\bibitem{L1} Lusztig, G.,  Finite dimensional Hopf algebras arising from quantized 
universal enveloping algebras, \emph{J. Amer. Math. Soc.} \textbf{3} (1990), 257--296.   

\bibitem{L2} Lusztig, G.,  Quantum groups  at  roots  of 1, \emph{ Geom. Dedicata} \textbf{35} (1990),  89--114.  

\bibitem{Lu} Lusztig, G.,  \emph{Introduction to quantum groups}, Birkh\"auser, 1993.

\bibitem{Ma3} Majid, S., \emph{Foundations of quantum group theory}, Cambridge Univ. Press, 1995.

\bibitem{M4} Masuoka, A., Abelian and non-abelian second cohomologies of quantized
enveloping algebras, \emph{J.  Algebra} \textbf{320} (2008), 1--47.


\bibitem{MS} Milinski, A.,  Schneider, H.-J.,   Pointed indecomposable Hopf algebras over Coxeter groups, 
\emph{Contemp. Math.} \textbf{267} (2000), 215--236.
  
\bibitem{Momb} Mombelli, M., Families of finite-dimensional Hopf algebras with the Chevalley property,
\emph{Algebr. Represent. Theory} \textbf{16} (2013), 421--435. 

  
\bibitem{Mo1} Montgomery, S., \emph{Hopf algebras and their actions on rings}, CMBS \textbf{82},  Amer. Math. Soc. 1993.

\bibitem{MW} Montgomery, S., Whiterspoon, S., Irreducible representations of
crossed products, \emph{J. Pure Appl. Algebra} \textbf{129} (1998), 315--326.


\bibitem{N4-mem} Natale, S.,  \emph{Semisolvability of semisimple Hopf algebras of low dimension},
Mem. Amer. Math. Soc. \textbf{186}, 123 pp. (2007).


\bibitem{Nat-ext} Natale, S.,  Hopf algebra extensions of group algebras and Tambara-Yamagami categories, 
\emph{Algebr. Represent. Theory} \textbf{13} (2010), 673--691.

\bibitem{N4-60} Natale, S.,  Semisimple Hopf algebras of dimension 60, \emph{J. Algebra} \textbf{324} (2010),
3017--3034.

\bibitem{Ng} Ng, S.-H., Non-semisimple Hopf algebras of dimension $p\sp 2$,  \emph{J. Algebra} \textbf{255} (2002),  182--197.

\bibitem{N}  Nichols, W. D., Bialgebras of type one, \emph{Comm. Algebra} \textbf{6} (1978),  1521--1552.  

\bibitem{Nik-tw}  Nikshych, D., $K_0$-rings and twisting of finite-dimensional semisimple
Hopf algebras, \emph{Comm. Algebra} \textbf{26} (1998), 321--342; Corrigendum, \textbf{26} (1998), 2019.

\bibitem{Nik}  Nikshych, D., Non-group-theoretical semisimple Hopf algebras from group actions on fusion categories. 
\emph{Selecta Math. (N. S.)} \textbf{14} (2008), 145--161.

\bibitem{Rad} Radford, D. E., \emph{Hopf algebras}, Series on Knots and Everything 49. Hackensack, NJ: World Scientific. xxii, 559 p.  (2012).

\bibitem{Ro2} Rosso, M., Quantum groups and quantum shuffles, \emph{Invent. Math.} \textbf{133}   (1998),  399--416.

\bibitem{S-bigal} Schauenburg, P., Hopf bi-Galois extensions, \emph{Comm. Algebra} \textbf{24} (1996), 3797--3825.

\bibitem{Sbg} Schauenburg, P., A characterization of the Borel-like subalgebras of quantum enveloping algebras, \emph{Comm. Algebra} {\bf 24} (1996),  2811--2823.

\bibitem{Sch2}  Schneider, H.-J., {\it  Lectures on Hopf algebras},  Trab. Mat. 31/95 (FaMAF), 1995.


\bibitem{Ser-superKM} Serganova, V., Kac-Moody superalgebras and integrability, \emph{Progr. Math.} \textbf{288} (2011), 169--218.

\bibitem{St-finite-ss} \c{S}tefan, D., The set of types of n-dimensional semisimple and cosemisimple Hopf algebras is finite, 
\emph{J. Algebra} \textbf{193} (1997), 571--580.

\bibitem{St} \c{S}tefan, D., Hopf algebras of low dimension, \emph{J. Algebra} \textbf{211} (1999), 343--361.

\bibitem{Sw}  Sweedler, M. E., \emph{Hopf algebras}, Benjamin, New York, 1969.

\bibitem{V} Vendramin, L., Nichols algebras associated to the transpositions of the symmetric group are twist-equivalent, 
\emph{Proc. Amer. Math. Soc.} \textbf{140} (2012),  3715--3723.

\bibitem{W} Woronowicz, S. L., Differential calculus on compact matrix pseudogroups (quantum
  groups), \emph{Comm. Math. Phys.} \textbf{122} (1989),  125--170.

\bibitem{Z} Zhu,  Y., Hopf algebras  of prime
dimension, \emph{Intern. Math. Res. Not.} {\bf 1} (1994),  53--59.
 
\end{thebibliography}
\end{document}